\documentclass{article}

\usepackage{latexsym, amssymb}

\newcommand\newsection[1]{\bigskip\bigskip\refstepcounter{section}
\noindent{\large\bf\thesection\ #1}}
\newcommand\newsubsection[1]{\medskip\refstepcounter{subsection}
{\noindent\bf\thesubsection\ #1.\ }}
\newcommand\newfigur[2]{\medskip\center\refstepcounter{subsection}
{\bf #1 \thesubsection.\ } #2.\ }

\newenvironment{theorem}{\newsubsection{Theorem}\sl}{}
\newenvironment{corollary}{\newsubsection{Corollary}\sl}{}
\newenvironment{lemma}{\newsubsection{Lemma}\sl}{}
\newenvironment{example}{\newsubsection{Example}}{}
\newenvironment{figur}[1]{\newfigur{Figure}{#1}}{}
\newenvironment{definition}{\newsubsection{Definition}}{}

\newenvironment{remark}{\newsubsection{Remark}}{}
\newenvironment{proof}{\medskip\noindent{\it Proof.\ }}{\mbox{$\Box$}}

\newcommand\intv[3]{{#1}\leq{#2}\leq{#3}}
\newcommand\interval[3]{{#1}\leq{#2}\leq{#3}}

\newcommand\be{\begin{equation}}
\newcommand\ee{\end{equation}}
\newcommand\hs[1]{\hspace{#1 em}}
\newcommand\bea{\begin{eqnarray}}
\newcommand\eea{\end{eqnarray}}

\def\reals{{\mathbb{R}}}
\def\complexes{{\mathbb{C}}}

\begin{document}

\centerline{\bf The Terwilliger Algebra of an Almost-Bipartite}
\centerline{\bf P- and Q-polynomial Association Scheme}

\bigskip

\centerline{John S. Caughman, Mark S. MacLean, and Paul M. Terwilliger}

\bigskip

\bigskip


\noindent{\bf Abstract.\ }
Let $Y$ denote a $D$-class symmetric association
scheme with $D \geq 3$, and suppose $Y$ is almost-bipartite
P- and Q-polynomial.  Let $x$ denote a vertex of $Y$ and let 
$T=T(x)$ denote the corresponding Terwilliger algebra.
We prove that any irreducible $T$-module $W$
is both thin and dual thin in the sense of Terwilliger.  We produce
two bases for $W$ and describe the action of $T$ on these bases.
We prove that the isomorphism class of $W$  as a $T$-module
is determined by two parameters, the dual endpoint and diameter of $W$.
We find a recurrence which gives the multiplicities with which
the irreducible $T$-modules occur in the standard module.
We compute this multiplicity for those
irreducible $T$-modules which have diameter at least $D-3$.

\bigskip



\newsection{Introduction}

The Terwilliger algebra of a commutative association scheme was introduced
    in \cite{terwSubconstituentI}. This algebra is a finite-dimensional,
    semisimple $\complexes$-algebra, and is noncommutative in general.
    The Terwilliger algebra has been used to study P- and Q-polynomial
    schemes \cite{caughmanModules}, \cite{terwSubconstituentI}, 
    group schemes \cite{balmacedaOura},
    \cite{bannaiMunemasa}, strongly regular graphs \cite{TomiyamaYamazaki},
    Doob schemes \cite{tanabe}, and schemes
    over the Galois rings of characteristic four \cite{ishibashi}.
    Other work involving this algebra can be found in
\cite{collins},
    \cite{collinsgirth},
    \cite{curtinBipartite}, \cite{curtinTerw},
    \cite{dickieThin}, \cite{egge}, \cite{go}, \cite{gotight},
    \cite{hobartIto}, \cite{ishi-ito-yamada}, 
    \cite{terwSubconstituentII},
    \cite{terwSubconstituentIII},
    and \cite{watatani}.

    The Terwilliger algebra is particularly well-suited for
    studying P- and Q-polynomial schemes; nevertheless, it is 
    apparent from \cite{terwSubconstituentI}
    that the intersection numbers of these schemes do not
    completely determine the structure of the algebra.
    In this article, we consider the Terwilliger algebra of
    an almost-bipartite P- and Q-polynomial scheme. We show that with the
    added almost-bipartite assumption, the intersection numbers of the
    scheme completely determine the structure of the algebra.

    To describe our results, let $Y=(X, \{ R_i \}_{\interval0iD})$
    denote a symmetric association scheme with $D \geq 3$.
    Suppose $Y$ is almost-bipartite P- and Q-polynomial.
    Fix any $x \in X$, and let $T=T(x)$ denote
    the Terwilliger algebra of $Y$ with respect to $x$.
    $T$ acts faithfully on the vector space $\complexes^{X}$ by
    matrix multiplication; we refer to $\complexes^{X}$ as the
    {\em standard module}. Since $T$ is semisimple,
    $\complexes^{X}$ decomposes into a
    direct sum of irreducible $T$-modules.

    Let $W$ denote an irreducible $T$-module contained in
    $\complexes^{X}$.  We show that $W$ is thin and dual thin in
    the sense of Terwilliger (Lemma \ref{his3.8new}).
    We produce two bases for $W$ with respect to which
    the action of $T$ is particularly simple (Theorem \ref{basis}).
    To describe this action, we use two sets of scalars, the {\em
    intersection numbers} of $W$ and the {\em dual intersection numbers}
    of $W$.  We compute these scalars in terms of the
    eigenvalues of $Y$, the dual eigenvalues of $Y$, and two additional
    parameters, called the {\em dual endpoint} and {\em diameter} of $W$
    (Theorems \ref{cwbw-theta}, \ref{cswbsw-theta}).
    We show that the dual endpoint and diameter of $W$
    determine its isomorphism class as a $T$-module (Theorem \ref{isom}).

    Combining our above results, we find a recurrence which
    gives the multiplicities with which the irreducible $T$-modules
    occur in $\complexes^{X}$ (Theorem \ref{bigkahoona}).  
    We compute this multiplicity for those
irreducible $T$-modules which have diameter at least $D-3$.   (Example
    \ref{multlisting}).

        In a future paper, we intend to use these results to
    study the subconstituents of an almost-bipartite P- and Q-polynomial
    scheme.  We hope this will produce a classification of
    these schemes.  Our work is closely related to that
    of B. Curtin concerning $2$-thin distance-regular graphs
    \cite{curtinBipartite}.


\newsection{Association Schemes}

\begin{definition}\label{schemedef}
    By a {\em symmetric association scheme} (or {\em scheme} for short) we
    mean a pair
    $Y=(X, \{ R_i \}_{\interval0iD})$, where $X$ is a nonempty finite
    set, $D$ is a nonnegative integer, and $R_0$, $\ldots$, $R_D$ are
    nonempty subsets of $X \times X$ such that
    \begin{enumerate}
    \begin{item}\label{tallyhoo}
    $ \{ R_i \} _{\interval0iD} $ is a partition of $X \times X$;
    \end{item}
    \begin{item}
    $ R_0 = \{ xx\; |\; x \in X \} $;
    \end{item}
    \begin{item}
    $ R_i^t = R_i$ for $\interval 0iD$,
    where $R_i^t =\{ yx\; |\; xy \in R_i \}$;
    \end{item}
    \begin{item}\label{tallyho}
    For all $h,i,j$ $(\intv0{h,i,j}D)$, and for all $x,y \in X$ such that
    $xy \in R_h$, the scalar
$$p^h_{ij} = | \{ z \in X | xz \in R_i, \mbox{ and } zy \in R_j \} | $$
    is independent of $x,y$.
    \end{item}
    \end{enumerate}
    The constants $p^h_{ij}$ are called the {\em intersection
    numbers of $Y$}.

\end{definition}

    For the rest of this section, let
    $Y=(X, \{ R_i \}_{\interval0iD})$ denote a scheme.  We begin with
    a few comments about the intersection numbers of $Y$.
    For all integers $i$ $(\intv0iD)$, set $k_i := p^0_{ii}$,
    and note that $k_i \not= 0$ since $R_i$ is nonempty.
    We refer to $k_i$ as the $i^{th}$ {\em valency} of $Y$.
    Observe that $p^0_{ij} = \delta_{ij} k_i$ $(\intv0{i,j}D)$.


We now recall the Bose-Mesner algebra of $Y$.
    Let Mat$_X(\complexes$) denote the $\complexes$-algebra of
    matrices with entries in $\complexes$, where the rows and columns
    are indexed by $X$.
    For each integer $i \; (\interval0iD)$, let
    $A_i$ denote the matrix in Mat$_X(\complexes$) with
    $xy$-entry
    \begin{equation}\label{assdef}
    (A_i)_{xy} = \left\{
    \begin{array}{lr} 1 & \mbox{if } xy \in R_i, \\
                      0 & \mbox{if } xy \not\in R_i
    \end{array} \right. \; \; \; \; \; \; \;  (x,y \in X).
    \end{equation}
    We refer to $A_i$ as the  $i^{th}$ {\em associate matrix of $Y$}.
    By Definition \ref{schemedef}, the associate matrices satisfy:
    (i) $A_0 = I$, where $I$ is the identity matrix in Mat$_X(\complexes$);
    (ii) the conjugate-transpose $\overline{A}_i^t = A_i$ $(\interval0iD)$;
    (iii) $A_0 + A_1 + \cdots + A_D = J$,
    where $J$ is the all 1's matrix in Mat$_X(\complexes$);
    (iv) $A_i A_j = \sum_{h=0}^{D} p^h_{ij} A_h$ $(\interval0{i,j}D)$.

    It follows from (i)-(iv) that $A_0, ..., A_D$ form a basis for
    a subalgebra $M$ of Mat$_X(\complexes)$.  $M$ is known as the
    {\em Bose-Mesner algebra of $Y$}.  Observe that $M$ is commutative,
    since the associate matrices are symmetric.

    By \cite[p.45]{bcn}, the algebra $M$ has
    a second basis $E_0, ..., E_D$ satisfying:
    (i) $E_0 = |X|^{-1} J$;
    (ii) $\overline{E}_i^t = E_i$ $(\interval0iD)$;
    (iii) $E_i E_j = \delta_{ij} E_i$ $(\interval0{i,j}D)$;
    (iv) $E_0 + E_1 + \cdots + E_D = I$.
    We refer to $E_i$ as the  $i^{th}$ {\em primitive idempotent of $Y$}
    for $\interval0iD$. For convenience
    we define $E_i :=0$ for
    $i > D$ and $i<0$.

    For all integers $i$ $(\intv0iD)$, set $m_i := rank(E_i)$,
    and note that $m_i \not= 0$.
    We refer to $m_i$ as the $i^{th}$ {\em multiplicity of $Y$}.

    Since $A_0, ..., A_D$  and $E_0, ..., E_D$ are both bases for $M$,
    there exist complex scalars $p_i(j)$, $q_i(j)$ $(\interval0{i,j}D)$ which
    satisfy
\begin{eqnarray}
    A_i &=& \sum_{j=0}^D p_i(j) E_j \hspace{3em} (\interval0iD), \label{ww1}\\
    E_i &=& |X|^{-1} \sum_{j=0}^D q_i(j) A_j \hspace{3em}
                        (\interval0iD). \label{ww2}
\end{eqnarray}
    By \cite[pp. 59, 63]{bannai}, the $p_i(j)$, $q_i(j)$ are real.
    We refer to $p_i(j)$ (resp. $q_i(j)$) as the $j^{th}$ {\em eigenvalue}
    (resp. $j^{th}$ {\em dual eigenvalue}) associated with $A_i$
    (resp. $E_i$).
    By \cite[p.63]{bannai}, the eigenvalues and dual eigenvalues
    satisfy
\begin{displaymath}
    \frac{p_i(j)}{k_i} = \frac{q_j(i)}{m_j} \hspace{3em} (\intv0{i,j}D).
\end{displaymath}
    We now recall the Krein parameters of $Y$. Observe that
    $A_i \circ A_j = \delta_{ij} A_i$ $(\interval0{i,j}D)$,
    where $\circ$ denotes the entry-wise matrix product.  It follows that
    $M$ is closed under $\circ$, so there exist complex scalars
    $q^h_{ij}$ satisfying
    $E_i \circ E_j = |X|^{-1} \sum_{h=0}^{D} q^h_{ij} E_h$
    $(\interval0{i,j}D)$.
    The constants $q^h_{ij}$ are called the {\em Krein parameters}
    of $Y$.
    By \cite[pp.67-69]{bannai}, the Krein parameters are real, and
    $q^0_{ij} = \delta_{ij} m_i$ $(\intv0{i,j}D)$.


We now recall the dual Bose-Mesner algebra of $Y$.
    For the rest of this section, fix any $x \in X$.  For each
    integer $i \; (\interval0iD)$, let
    $E^*_i = E^*_i(x)$ denote the diagonal matrix in Mat$_X(\complexes$)
    with $yy$-entry
    \begin{equation}\label{e*def}
    (E^*_i)_{yy} = \left\{
    \begin{array}{lr} 1 & \mbox{if } xy \in R_i, \\
                      0 & \mbox{if } xy \not\in R_i
    \end{array} \right. \; \; \; \; \; \; \; (y \in X).
    \end{equation}
    We refer to $E^*_i$ as the $i^{th}$
    {\em dual idempotent of $Y$ with respect to x}.
    For convenience, set $E^*_i :=0$ if $i > D$ or $i<0$.
    From the definition, the dual idempotents satisfy:
    (i) $\overline{E^*_i}^{t} = E^*_i$ $(\interval0iD)$;
    (ii) $E^*_i E^*_j = \delta_{ij} E^*_i$ $(\interval0{i,j}D)$;
    (iii) $E^*_0 + E^*_1 + \cdots + E^*_D = I$.

    It follows that the matrices $E^*_0, ..., E^*_D$
    form a basis for a subalgebra
    $M^*=M^*(x)$ of Mat$_X(\complexes$). $M^*$ is known as
    the {\em dual Bose-Mesner algebra of} $Y$ {\em with respect to $x$.}
    Observe that $M^*$ is commutative since the dual
    idempotents are diagonal.

    For each integer $i \; (\interval0iD)$, let $A^*_i = A^*_i(x)$
    denote the diagonal matrix in Mat$_X(\complexes$) with $yy$-entry
    \begin{equation}\label{a*def}
    (A^*_i)_{yy} = |X| (E_i)_{xy}  \; \; \; \; (y \in X).
    \end{equation}
    We refer to $A_i^*$ as the $i^{th}$ {\em dual associate
    matrix of $Y$ with respect to x}.
    Combining (\ref{ww1}),(\ref{ww2}) with
    (\ref{e*def}), (\ref{a*def}),
\begin{eqnarray}
    A^*_i &=& \sum_{j=0}^D q_i(j) E^*_j \hspace{3em} (\interval0iD),
                                    \label{ww3}\\
    E^*_i &=& |X|^{-1} \sum_{j=0}^D p_i(j) A^*_j \hspace{3em}
                               (\interval0iD). \label{ww4}
\end{eqnarray}
    It follows that $A^*_0, ..., A^*_D$ form
    a second basis for $M^*$.

    From the definitions, the dual associate matrices satisfy:
    (i) $A^*_0 = I$;
    (ii) $\overline{A^*_i}^{t} = A^*_i$ $(\interval0iD)$;
    (iii) $A^*_i A^*_j = \sum_{h=0}^{D} q^h_{ij} A^*_h$ $(\interval0{i,j}D)$;
    (iv) $A^*_0 + A^*_1 + \cdots + A^*_D = |X| E^*_0$.
%
%

\newsection{The Terwilliger Algebra and its Modules}

    Let $Y=(X, \{ R_i \}_{\interval0iD})$ denote a scheme.
    Fix any $x \in X$, and write $M^*=M^*(x)$.
    Let $T=T(x)$ denote the subalgebra of Mat$_X(\complexes$)
    generated by $M$ and $M^*$.  We call $T$ the {\em Terwilliger
    algebra of $Y$ with respect to $x$}.

    In \cite[Lemma 3.2]{terwSubconstituentI}, it is shown that
    for all integers $h$,$i$,$j \; (\interval0{h,i,j}D)$,
\begin{eqnarray}
    p^h_{ij} =0 \hspace{2em} &\mbox{if and only
                    if}& \hspace{2em} E^*_i A_j E^*_h =0,
    \label{triangle} \\
    q^h_{ij} =0 \hspace{2em} &\mbox{if and only
                    if}& \hspace{2em} E_i A^*_j E_h =0,
    \label{triangleQ}
\end{eqnarray}
    where $A^*_i = A^*_i(x)$, $E^*_i = E^*_i(x) \;
    (\interval0iD)$.

    Let $V$ denote the vector
    space $\complexes ^{X}$ (column vectors),
    where the coordinates are indexed by $X$.  Then
    Mat$_X(\complexes$) acts on $V$ by left multiplication.
    We endow $V$ with the inner product $\langle \; , \; \rangle$  satisfying
    $\langle u , v \rangle := u^t \overline{v}$ for all $u, v \in V$.
    Observe that $V = \sum_{i=0}^D E_i V$ (orthogonal direct sum).
    Similarly, we have the decomposition
    $V = \sum_{i=0}^D E^*_i V$ (orthogonal direct sum).

    By a {\em $T$-module}, we mean a subspace $W$ of $V$ such that
    $TW \subseteq W$. We refer to $V$ itself as the
    {\em standard module} for $T$.
    Let $W$, $W'$ denote $T$-modules.
    By a {\em $T$-module isomorphism} from $W$ to
    $W'$, we mean an isomorphism of vector
    spaces $\phi: W \rightarrow W'$ such that
\begin{displaymath}
    (B \phi - \phi B)W = 0 \hspace{3em} (\forall B \in T).
\end{displaymath}
    $W$, $W'$ are said to be {\em $T$-isomorphic} whenever there exists
    a $T$-module isomorphism from $W$ to $W'$.
    A $T$-module $W$ is said to be {\em irreducible}
    whenever $W \not= 0$ and $W$ contains
    no $T$-modules other than $0$ and $W$.

    Because $T$ is closed under the conjugate-transpose map, $T$ is
    semisimple.  It follows that for any $T$-module $W$
    and any $T$-module $U \subseteq W$ there exists a unique
    $T$-module $U' \subseteq W$ such that
\begin{displaymath}
    W = U + U' \hspace{3em} (\mbox{orthogonal direct sum}).
\end{displaymath}
    Moreover, $W$ is an orthogonal
    direct sum of irreducible $T$-modules.

    Now let $W$ denote an irreducible $T$-module.  Observe that
\begin{equation}\label{WE*decomp}
    W = \sum E^*_i W \hspace{3em} \mbox{(orthogonal direct sum),}
\end{equation}
    where the sum is taken over all the indices $i$ $(\intv0iD)$ such
    that $E^*_iW \not=0$.
    We set $d := | \{ i \; | \; E^*_i W \not= 0 \} | -1$,
    and observe that the dimension of $W$ is at least $d+1$.
    We refer to $d$ as the {\em diameter} of $W$.
    $W$ is said to be {\em thin} whenever
    $dim (E^*_i W) \leq 1$ $(\interval0iD)$.
    Note that $W$ is thin if and only if the diameter of $W$
    equals $dim(W) -1$.

Similarly,
\begin{displaymath}
    W = \sum E_i W \hspace{3em} \mbox{(orthogonal direct sum),}
\end{displaymath}
    where the sum is taken over all the indices $i$ $(\intv0iD)$ such
    that $E_iW \not=0$.  We set
    $d^* := | \{ i \; | \; E_i W \not= 0 \} | -1$,
    and observe that the dimension of $W$ is at least $d^*+1$.
    We refer to $d^*$ as the {\em dual diameter} of $W$.
    $W$ is said to be {\em dual thin} whenever
    $dim (E_i W) \leq 1$ $(\interval0iD)$.
    Note that $W$ is dual thin if and only if the dual diameter of
    $W$ equals $dim(W)-1$.
    
    We wish to emphasize the following point, which follows 
    immediately from the above discussion.
    
\begin{lemma}  \label{diams=}
    Let $Y=(X, \{ R_i \}_{\interval0iD})$ denote a scheme.  Fix any 
    $x \in X$, and write $T=T(x)$.  Let $W$ denote an irreducible 
    $T$-module that is both thin and dual thin.  
    Then the diameter and dual diameter of $W$ are equal.
\end{lemma}

%
%

\newsection{The P-Polynomial Property}

    Let $Y=(X, \{ R_i \}_{\interval0iD})$ denote a scheme.
    We say that $Y$ is {\em P-polynomial} (with
    respect to the ordering $R_0, ..., R_D$ of the associate classes)
    whenever for all integers $h,i,j \; (\interval0{h,i,j}D)$,
    \begin{eqnarray}
    &p^h_{ij} = 0 \mbox{
    {\sl if one of $h,i,j$ is greater than the sum of
    the other two,}}& \label{ppolydef1} \\
    &p^h_{ij} \not=  0 \mbox{
    {\sl if one of $h,i,j$ equals the sum of
    the other two.}}&
    \end{eqnarray}
    For the rest of this section, assume $Y$ is P-polynomial.
    We abbreviate $c_i:=p^i_{1i-1}$ $(\interval1iD)$,
    $a_i:=p^i_{1i}$ $(\interval0iD)$,
    $b_i:=p^i_{1i+1}$ $(\interval0i{D-1})$,
    and define $c_0 := 0, \, b_D := 0.$ We note that $a_0=0$ and $c_1=1$.
    By \cite[Prop. III.1.2]{bannai},
\begin{displaymath}
    k_i = \frac{b_0 b_1 \cdots b_{i-1}}
        {c_1 c_2 \cdots c_i} \hspace{3em} (\interval0iD).
\end{displaymath}
    Of particular interest are the matrix $A := A_1$ and the scalars
    $\theta_i:=p_1(i)$ $(\interval0iD)$.
    It follows from (\ref{ww1}) that
    \begin{equation}\label{eieigsp}
    A E_i = \theta_i E_i \hspace{3em} (\interval0iD).
    \end{equation}
    It is shown in \cite[p.190]{bannai} that
    \begin{equation}\label{thdist}
    \theta_i \not= \theta_j \hspace{2em} \mbox{if} \hspace{2em} i \not= j
    \hspace{3em} (\interval0{i,j}D),
    \end{equation}
    and also that $A_i = v_i (A)$ $(\intv0iD)$,
    where $v_i$ is a polynomial with real coefficients and degree
    exactly $i$.
    In particular, $A$ multiplicatively generates $M$, the Bose-Mesner algebra.
    By (\ref{ww1}), it follows that
\begin{displaymath}
    p_i(j) = v_i (\theta_j) \hspace{3em} (\intv0{i,j}D).
\end{displaymath}


We now recall the raising, lowering, and flat matrices of
    $Y$.  Fix any $x \in X$, and write $E^*_i = E^*_i(x)$ $(\interval0iD)$.
    Define matrices $R=R(x), \; F= F(x), \; L=L(x)$ by
\begin{equation}\label{fdef}
    R := \sum^{D}_{i=0} E^*_{i+1} A E^*_i, \hspace{2em}
    F := \sum^{D}_{i=0} E^*_i A E^*_i, \hspace{2em}
    L := \sum^{D}_{i=0} E^*_{i-1} A E^*_i.
\end{equation}
    Note that $R$, $F$, and $L$ have real entries by
    (\ref{assdef}),(\ref{e*def}).  Also, observe that $F$ is
    symmetric and $R= L^t$.
    By (\ref{triangle}) and (\ref{ppolydef1}),
\begin{equation}\label{sumrfl}
A = R + F + L.
\end{equation}
Using (\ref{fdef}) and recalling $E^{*}_{-1}=0$, $E^{*}_{D+1}=0$, we find
\be  \label{hick}
  RE^*_i = E^*_{i+1}R \quad (-1 \leq i \leq D), \qquad  FE^*_i = E^*_{i}F 
  \quad (\intv0iD), \qquad  LE^*_i = E^*_{i-1}L \quad (\intv0i{D+1}).
  \ee


\newsection{The $T$-Modules of P-Polynomial Schemes}

    In this section, we describe the
    irreducible $T$-modules of P-polynomial schemes.

    Let $Y=(X, \{ R_i \}_{\interval0iD})$ denote a scheme which
    is P-polynomial with respect to the ordering
    $R_0,...,R_D$ of the associate classes.  Fix
    any $x \in X$ and write $T=T(x)$.  Let $W$ denote an
    irreducible $T$-module.  We define the {\em endpoint} $r$ of
    $W$ by
\begin{eqnarray*}
    r &:=& min \{ i \; | \; 0 \leq i \leq D, \; E^*_i W \not= 0 \}. 
\end{eqnarray*}
    We observe that $0 \leq r \leq D-d, \;$ where $d$ denotes the
    diameter of $W$.

    In \cite[Lemma 3.9]{terwSubconstituentI}, it was shown that
    $R E^*_i W \not= 0$ $(r \leq i < r+d)$, $L E^*_{i} W \not= 0$
    $(r < i \leq r+d)$,
    and also that
\begin{eqnarray}\label{ddefeq}
    E^*_i W &\not= & 0  \hspace{2em} \mbox{ iff } \hspace{2em}
    \interval{r}i{r+d} \hspace{3em} (\interval0iD).
\end{eqnarray}
    By \cite[Lemma 5.1]{caughmanModules}, we have
\begin{eqnarray*}
    2r+d^* \geq D,
\end{eqnarray*}
    where $d^{*}$ denotes the dual diameter of $W$.
    
    \bigskip

\begin{lemma}\label{ppolythin}
    \cite[Lemma 3.9]{terwSubconstituentI}
    Let $Y=(X, \{ R_i \}_{\interval0iD})$ denote a scheme which
    is P-polynomial with respect to the ordering $R_0,...,R_D$
    of the associate classes.  Fix any $x \in X$,
    and write $E^*_i = E^*_i(x) \; (\interval0iD)$,
    $T=T(x)$.  Let $W$ denote a thin, irreducible $T$-module
    with endpoint $r$. Then
\begin{enumerate}
\begin{item} $W = M E^*_r W.$
\end{item}
\begin{item} $E_i W = E_i E^*_r W \hspace{3em} (\interval0iD).$
\end{item}
\begin{item} $W$ is dual thin.
\end{item}
\end{enumerate}
\end{lemma}


\newsection{The Q-Polynomial Property}

    Let $Y=(X, \{ R_i \}_{\interval0iD})$ denote a scheme.
    We say that $Y$ is {\em Q-polynomial} (with
    respect to the given ordering $E_{0}, E_{1}, \ldots, E_{D}$
    of the primitive idempotents)
    whenever for all integers $h,i,j \; (\interval0{h,i,j}D)$, 
    the Krein parameters satisfy 
    \begin{eqnarray*}
    &q^h_{ij} = 0 \mbox{
    {\sl if one of $h,i,j$ is greater than the sum of
    the other two,}}&  \\
    &q^h_{ij} \not=  0 \mbox{
    {\sl if one of $h,i,j$ equals the sum of
    the other two.}}&
    \end{eqnarray*}

    For the rest of this section, assume $Y$ is Q-polynomial with
    respect to the ordering $E_0,..., E_D$.
    We abbreviate $c^*_i:=q^i_{1i-1}$ $(\interval1iD)$,
    $a^*_i:=q^i_{1i}$ $(\interval0iD)$,
    $b^*_i:=q^i_{1i+1}$ $(\interval0i{D-1})$,
    and define $c^*_0 := 0, \, b^*_D := 0.$ We note that $a^*_0=0$ and  $c^*_1=1$
    \cite[Prop. II.3.7]{bannai}.
    By \cite[p.196]{bannai},
\begin{displaymath}
    m_i = \frac{b^*_0 b^*_1 \cdots b^*_{i-1}}{c^*_1 c^*_2 \cdots c^*_i}
         \hspace{3em} (\interval0iD).
\end{displaymath}
    Fix any $x \in X$ and write $E^*_i = E^*_i(x)$, $A^*_i = A^*_i(x)$
    $(\interval0iD)$.
    Of particular interest are the matrix
    $A^* := A^*_1(x)$ and the scalars
    $\theta^*_i := q_1(i)$ $(\interval0iD)$.
    By (\ref{ww3}),
    \begin{equation}\label{eieigsp*}
    A^* E^*_i = \theta^*_i E^*_i \hspace{3em} (\interval0iD).
    \end{equation}
    It is shown in \cite[p.193]{bannai} that
    \begin{equation}\label{thsdist}
    \theta^*_i \not= \theta^*_j \hspace{2em} \mbox{if}
    \hspace{2em} i \not= j \hspace{3em} (\interval0{i,j}D),
    \end{equation}
    and also that $A^*_i = v^*_i (A^*)$ $(\intv0iD)$,
    where $v^*_i$ is a polynomial with real coefficients and
    degree exactly $i$. In particular, $A^*$ generates
    the dual Bose-Mesner algebra $M^* = M^*(x)$.
    By (\ref{ww3}), it follows that
\begin{displaymath}
    q_i(j) = v^*_i (\theta^*_j) \hspace{3em} (\intv0{i,j}D).
\end{displaymath}


We now recall the dual raising, lowering, and flat matrices of $Y$.
    Define the matrices $R^*=R^*(x)$, $F^*= F^*(x), \; L^*=L^*(x)$ by
\begin{equation}\label{tuggy}
    R^* := \sum^{D}_{i=0} E_{i+1} A^* E_i, \hspace{2em}
    F^* := \sum^{D}_{i=0} E_i A^* E_i, \hspace{2em}
    L^* := \sum^{D}_{i=0} E_{i-1} A^* E_i.
\end{equation}
    Note that $R^*$, $F^*$, and $L^*$ have real entries by
    (\ref{ww3}), and since the $q_i(j)$ are real. Also, observe
    that $F^*$ is symmetric and $R^*= L^{*t}$.  Moreover,
\begin{equation}\label{sumrfl*}
A^* = R^* + F^* + L^*.
\end{equation}
Using (\ref{tuggy}) and recalling $E_{-1}=0$, $E_{D+1}=0$, we find
\be \label{hick-s}
  R^{*}E_i = E_{i+1}R^{*} \; (\intv{-1}iD), \quad F^{*}E_i = E_{i}F^{*}
  \;(\intv0iD), \quad L^{*}E_i = E_{i-1}L^{*} \;
  (\intv0i{D+1}).
  \ee


\newsection{The $T$-Modules of Q-Polynomial Schemes}

    In this section, we describe the
    irreducible $T$-modules of Q-polynomial schemes.

    Let $Y=(X, \{ R_i \}_{\interval0iD})$ denote a scheme which
    is Q-polynomial with respect to the ordering $E_0,...,E_D$
    of the primitive idempotents.  Fix any $x \in X$ and
    write $T=T(x)$.  Let $W$ denote an
    irreducible $T$-module.  We define the {\em dual endpoint } $ t$ of
    $W$ by
\begin{eqnarray*}
    t &:=& min \{ i \; | \; 0 \leq i \leq D, \; E_i W \not= 0 \}.
\end{eqnarray*}
    We observe that $0 \leq t \leq D-d^*, \;$ where $d^*$ denotes
    the dual diameter of $W$.

    In \cite[Lemma 3.12]{terwSubconstituentI}, it was shown that
    $R^{*} E_i W \not= 0$ $(t \leq i < t+d^{*})$, $L^{*} E_{i} W \not= 0$
    $(t < i \leq t+d^{*})$, 
    and also that
    \begin{eqnarray}
    E_i W &\not= & 0 \hspace{2em} \mbox{ iff } \hspace{2em}
    \interval{t}i{t+d^*}
    \hspace{3em} (\interval0iD). \label{d*defeq}
    \end{eqnarray}
    By \cite[Lemma 7.1]{caughmanModules}, we have 
\begin{eqnarray}\label{2t+deq}
    2t+d \geq D,
\end{eqnarray}
    where $d$ denotes the diameter of $W$.
    
    \bigskip

\begin{lemma}\label{qpolythin}\cite[Lemma 3.12]{terwSubconstituentI}
    Let $Y=(X, \{ R_i \}_{\interval0iD})$ denote a scheme which
    is Q-polynomial with respect to the ordering $E_0,...,E_D$
    of the primitive idempotents. Fix any $x \in X$, and write
    $E^*_i = E^*_i(x)$ $(\interval0iD)$, $M^* = M^*(x)$, $T=T(x)$.
    Let $W$ denote a dual thin, irreducible
    $T$-module with dual endpoint $t$.  Then (i)--(iii) hold below.
\begin{enumerate}
\begin{item} $W = M^* E_t W.$
\end{item}
\begin{item} $E^*_i W = E^*_i E_t W \hspace{3em} (\interval0iD).$
\end{item}
\begin{item} $W$ is thin.
\end{item}
\end{enumerate}
\end{lemma}


\newsection{The $T$-Modules of P- and Q-Polynomial Schemes}

    Let $Y=(X, \{ R_i \}_{\interval0iD})$ denote a scheme which
    is P-polynomial with respect to the ordering
    $R_0,...,R_D$ of the associate classes, and Q-polynomial 
    with respect to the 
    ordering $E_0,...,E_D$
    of the primitive idempotents.  Fix any $x \in X$ and write 
    $T=T(x)$.  Let $W$ denote a thin irreducible 
    $T$-module.  Observe $W$ is dual thin by Lemma \ref{ppolythin},
    and the diameter and dual diameter of $W$ coincide by Lemma 
    \ref{diams=}.  We now present two bases for $W$, one of which 
    diagonalizes $A$ and the other diagonalizes $A^{*}$.  We then 
    consider the action of $T$ on these bases.
    
%
%
%
%

\begin{theorem}\label{basis}
    Let $Y=(X, \{ R_i \}_{\interval0iD})$ denote a scheme which
    is P-polynomial with respect to the ordering
    $R_0,...,R_D$ of the associate classes, and
    Q-polynomial with respect to the ordering $E_0,...,E_D$
    of the primitive idempotents. Fix any $x \in X$, and write
    $E^*_i = E^*_i(x)$ $(\interval0iD)$, $T=T(x)$.  
    Let $W$ denote a thin irreducible $T$-module with
    endpoint $r$, dual endpoint $t$, and diameter $d$.
\begin{enumerate}
\begin{item} For all nonzero $v \in E_t W$, the vector
$E^*_i v$ is a basis for $E^*_i W$ for $\intv{r}{i}{r+d}$.
Moreover, $E^*_r v , E^*_{r+1} v, \ldots, E^*_{r+d}v$
    is a basis for $W$.
\end{item}
\begin{item} For all nonzero $v \in E^*_r W$, the vector
$E_i v$ is a basis for $E_i W$ for $\intv{t}{i}{t+d}$.
Moreover, $E_t v$,$E_{t+1} v$,$\ldots$, $E_{t+d}v$ is a basis
for $W$.
\end{item}
\end{enumerate}
\end{theorem}

\begin{proof}
(i).  Recall $W$ is dual thin by Lemma \ref{ppolythin} so $v$ spans 
$E_{t}W$.
    Fix any $i$ $(\intv{r}i{r+d})$, and observe
    $E^*_i W \not= 0$ by (\ref{ddefeq}).
    Also $E^*_i v$ spans $E^*_i W$, since by Lemma
    \ref{qpolythin} and the construction,
\begin{eqnarray}
    E^*_i W &=& E^*_i E_t W \nonumber \\
            &=& span(E^*_i v). \nonumber
\end{eqnarray}
    We have now shown that $E^*_i v$ is a basis
    for $E^*_i W$.  Applying (\ref{WE*decomp}), (\ref{ddefeq}),
    we find
    $E^*_r v,...,E^*_{r+d}v$ is a basis for $W$.

\noindent
(ii).  Similar to the proof of (i).
\end{proof}


\begin{definition}\label{intfunny8}
    Let $Y=(X, \{ R_i \}_{\interval0iD})$ denote a scheme which
    is P-polynomial with respect to the ordering
    $R_0,...,R_D$ of the associate classes, and Q-polynomial 
    with respect to the ordering $E_0,...,E_D$
    of the primitive idempotents. Fix any $x \in X$, and write
    $E^*_i = E^*_i(x)$ $(\interval0iD)$, $T=T(x)$. 
    Let $W$ denote a thin irreducible $T$-module with endpoint $r$, 
    dual endpoint $t$, and diameter $d$.  
    For all $i$ $(\intv0id)$, let $c_i(W)$, $a_{i}(W)$, $b_i(W)$ denote the
    complex scalars such that
    \begin{eqnarray}
        R E^*_{r+i-1} v  &=& c_i (W) E^*_{r+i} v, \label{cwiparteq} \\
        F E^{*}_{r+i} v &=& a_{i}(W) E^*_{r+i} v, \label{awiparteq} \\
        L E^*_{r+i+1} v  &=& b_i (W) E^*_{r+i} v,\label{bwiparteq}
    \end{eqnarray}
    where $v$ is any nonzero vector in $E_t W$.
    Since $E_{t}W$ has dimension 1, we see 
    $c_i(W)$, $a_{i}(W)$, $b_i(W)$ are independent of the choice of $v$.
    We refer to the $c_i(W)$, $a_{i}(W)$, $b_i(W)$ as the 
    {\em intersection numbers}
    of $W$.  We observe $c_{0}(W)=0$, $b_{d}(W)=0$.   
By the {\em intersection matrix} of $W$, we mean the tridiagonal matrix
\[ B(W) := \left(
    \begin{array}{cccccc}
    a_0(W) & b_0(W)  &         &            &           & {\bf 0} \\
    c_1(W) & a_1(W)  & b_1(W)  &            &           &         \\
           & c_2(W)  & \cdot   & \cdot      &           &         \\
           &         & \cdot   & \cdot      & \cdot     &         \\
           &         &         & \cdot      & \cdot     &  b_{d-1}(W) \\
   {\bf 0} &         &         &            &c_d(W) &   a_{d}(W)     \\
    \end{array}
\right) . \]
\end{definition}


\begin{definition}\label{int*}
    Let $Y=(X, \{ R_i \}_{\interval0iD})$ denote a scheme which
    is P-polynomial with respect to the ordering
    $R_0,...,R_D$ of the associate classes, and Q-polynomial 
    with respect to the ordering $E_0,...,E_D$
    of the primitive idempotents. Fix any $x \in X$, and write
    $E^*_i = E^*_i(x)$ $(\interval0iD)$, $T=T(x)$. 
    Let $W$ denote a thin irreducible $T$-module with endpoint $r$, 
    dual endpoint $t$, and  diameter $d$.
    For all $i$ $(\intv0id)$, let $c^*_i(W)$, $a^*_i(W)$, $b^*_i(W)$
    denote the complex scalars such that
    \begin{eqnarray}
        R^* E_{t+i-1} v  &=& c^*_i (W) E_{t+i} v, \label{cswiparteq} \\
        F^* E_{t+i} v  &=& a^*_i (W) E_{t+i} v, \label{aswiparteq} \\
        L^* E_{t+i+1} v  &=& b^*_i (W) E_{t+i} v, \label{bswiparteq}
    \end{eqnarray}
    where $v$ is any nonzero vector in $E^*_r W$.
    Since $E^{*}_{r}W$ has dimension 1, we see
    $c^*_i(W)$, $a^*_i(W)$, $b^*_i(W)$ are independent of
    the choice of $v$.
    We refer to the
    $c^*_i(W)$, $a^*_i(W)$, $b^*_i(W)$
    as the {\em dual intersection numbers} of $W$.  We observe
    $c^*_0(W)=0$, $b^*_d(W)=0$.   
By the {\em dual intersection matrix} of $W$, we mean the tridiagonal
matrix
\[ B^*(W) := \left(
    \begin{array}{cccccc}
    a^*_0(W) &   b^*_0(W) &          &        &           &  {\bf 0} \\
    c^*_1(W) &   a^*_1(W) & b^*_1(W) &        &           &          \\
             &   c^*_2(W) & \cdot    &  \cdot &           &          \\
             &            & \cdot    &  \cdot &  \cdot    &          \\
             &            &          &  \cdot &  \cdot    & b^*_{d-1}(W) \\
    {\bf 0}  &            &          &        & c^*_d(W)  & a^*_d(W) \\
    \end{array}
\right) . \]
\end{definition}

\begin{lemma}\label{intmatlem}
    Let $Y=(X, \{ R_i \}_{\interval0iD})$ denote a scheme which
    is P-polynomial with respect to the ordering
    $R_0,...,R_D$ of the associate classes, and Q-polynomial 
    with respect to the ordering $E_0,...,E_D$
    of the primitive idempotents. Fix any $x \in X$, and write
    $E^*_i = E^*_i(x)$ $(\interval0iD)$, $T=T(x)$. 
    Let $W$ denote a thin irreducible $T$-module with endpoint $r$,
    dual endpoint $t$, and diameter $d$.
\begin{enumerate}
\begin{item}
$B(W)$ is the matrix representing multiplication by $A$ with respect
    to the basis $E^*_r v$, $E^*_{r+1} v$, ..., $E^*_{r+d}v$, where
    $v$ is any nonzero vector in $E_t W$.
\end{item}
\begin{item}
    $Diag(\theta^*_r, \theta^*_{r+1}, ...,
    \theta^*_{r+d})$ is the matrix representing multiplication by 
    $A^{*}$ with respect to the basis $E^*_r v$, $E^*_{r+1} v$, ..., $E^*_{r+d}v$, where
    $v$ is any nonzero vector in $E_t W$. 
\end{item}
\begin{item}
$B^*(W)$ is the matrix representing multiplication by $A^*$ with respect
    to the basis $E_t v$, $E_{t+1} v$, ..., $E_{t+d}v$,
    where $v$ is any nonzero vector in $E^*_r W$.
\end{item}
\begin{item}
     $Diag(\theta_t, \theta_{t+1}, ...,
    \theta_{t+d})$ is the matrix representing multiplication by $A$ 
    with respect
    to the basis $E_t v$, $E_{t+1} v$, ..., $E_{t+d}v$,
    where $v$ is any nonzero vector in $E^*_r W$.
\end{item}
\end{enumerate}
\end{lemma}

\begin{proof}
(i).  Immediate from (\ref{sumrfl}) and Definition \ref{intfunny8}.

\noindent
(ii).  Immediate from (\ref{eieigsp*}).  

\noindent
(iii).  Similar to the proof of (i).

\noindent
(iv).  Similar to the proof of (ii).
\end{proof}

\begin{corollary}   \label{b(w)eigs}
     Let $Y=(X, \{ R_i \}_{\interval0iD})$ denote a scheme which
    is P-polynomial with respect to the ordering
    $R_0,...,R_D$ of the associate classes, and Q-polynomial 
    with respect to the ordering $E_0,...,E_D$
    of the primitive idempotents. Fix any $x \in X$, and write
    $T=T(x)$. 
    Let $W$ denote a thin irreducible $T$-module with endpoint $r$,
    dual endpoint $t$, and diameter $d$.
\begin{enumerate}
\begin{item}  \label{bweigvals2}
The eigenvalues of $B(W)$ are $\theta_t$, $\theta_{t+1}$, ...,
    $\theta_{t+d}$.
\end{item}

\begin{item}  \label{bwseigvals2}
The eigenvalues of $B^*(W)$ are $\theta^*_r$, $\theta^*_{r+1}$, ...,
    $\theta^*_{r+d}$.
\end{item}
\end{enumerate}
\end{corollary}

\begin{corollary}\label{tracelem2}
    Let $Y=(X, \{ R_i \}_{\interval0iD})$ denote a scheme which
    is P-polynomial with respect to the ordering
    $R_0,...,R_D$ of the associate classes, and Q-polynomial 
    with respect to the ordering $E_0,...,E_D$
    of the primitive idempotents. Fix any $x \in X$, and write
    $E^*_i = E^*_i(x)$ $(\interval0iD)$, $T=T(x)$. 
    Let $W$ denote a thin irreducible $T$-module with endpoint $r$, dual endpoint $t$, 
    and diameter $d$.  Then (i), (ii) hold below. 
\begin{enumerate}
\begin{item} 
    $\displaystyle{\sum_{i=0}^d a_i(W) = \sum_{i=t}^{t+d} \theta_i.}$
\end{item}

\begin{item}
    $\displaystyle{\sum_{i=0}^d a_i^*(W) = \sum_{i=r}^{r+d} 
    \theta^*_i.}$
\end{item}
\end{enumerate}
\end{corollary}

\begin{proof}  (i).  By Corollary \ref{b(w)eigs}(\ref{bweigvals2}), 
both sides of the 
equation in (i)
    equal the trace of $B(W)$.
    
    \noindent
    (ii).   By Corollary \ref{b(w)eigs}(\ref{bwseigvals2}), 
    both sides of the equation 
    in (ii)
    equal the trace of $B^*(W)$.
\end{proof}


\newsection{Almost-Bipartite P- and Q-Polynomial Schemes}

    Let $Y=(X, \{ R_i \}_{\interval0iD})$ denote a scheme which
    is P-polynomial with respect to the ordering
    $R_0,...,R_D$ of the associate classes.  We say
    $Y$ is {\em almost-bipartite} (with respect to the 
    P-polynomial
    ordering) whenever $a_i = 0$ for $\intv0i{D-1}$ and $a_{D} \not= 
    0$.  
    
    For the remainder of this article, we shall be concerned
    with  P- and Q-polynomial schemes for which the 
    P-polynomial structure is almost-bipartite.  We thus make the 
    following definition.

\medskip


\begin{definition} \label{setup}
    Let $Y=(X, \{ R_i \}_{\interval0iD})$ denote a scheme with $D \geq 3$
    which is almost-bipartite P-polynomial with respect to the ordering
    $R_0, ... , R_D$ of the associate classes,
    and Q-polynomial with respect to the
    ordering $E_0,...,E_D$ of the primitive idempotents. Fix any $x \in X$,
    and write $T = T(x)$ to denote the Terwilliger algebra of $Y$ with
    respect to $x$.  (Where the context allows, we will also suppress
    the reference to $x$ for the individual matrices in $T$
     -- e.g., $E^*_0 = E^*_0(x)$, $R=R(x)$, etc.).
\end{definition}

\medskip


\begin{lemma}\label{flat}
Let $Y$ be as in Definition \ref{setup}. Then (i)--(iii) hold below.
\begin{enumerate}
\begin{item}
     $E_{i}^{*}AE_{i}^{*}=0$ $\;\;(\intv0i{D-1})$, and $E_{D}^{*}A E_{D}^{*} 
     \not= 0.$
\end{item}
\begin{item}
 $F = 
E_{D}^{*}A E_{D}^{*}$, where $F$ is the matrix from (\ref{fdef}).
\end{item}
\begin{item}
   $FE_{i}^{*}=0$ $\;\;(\intv0i{D-1})$.
\end{item}
\end{enumerate}
\end{lemma}

\begin{proof} (i).  Immediate from (\ref{triangle}) and the fact that 
    $a_{i} =0 \; \; (\intv0i{D-1})$, $a_{D} \not= 0$. 
    
  \noindent  (ii), (iii).  Immediate using (i).    
\end{proof}

\begin{theorem}  \cite[Theorem 14.3]{collins}  \label{width} 
    With reference to Definition \ref{setup}, let $W$ denote an 
    irreducible $T$-module with endpoint $r$ and diameter $d$.  Then 
    $r+d=D$.
\end{theorem}


\newsection{Each Irreducible $T$-Module is Thin and Dual Thin}
     \label{basessec}

    Let $Y$ be as in Definition \ref{setup}, and let $W$ denote
    an irreducible $T$-module.
    In this section, we show that $W$ is both thin and dual thin.
    

\begin{lemma}\label{1strel}
    With reference to Definition \ref{setup},
    let $W$ denote an irreducible $T$-module with dual endpoint $t$, and fix
    any nonzero $v \in E_{t} W$. Then
\begin{enumerate}
\begin{item}
    $\displaystyle{R E^*_{i-1} v + L E^*_{i+1} v = \theta_{t} E^*_i v
    \hspace{1em} (\interval0i{D-1}),}$
\end{item}
\begin{item} 
    $\displaystyle{RE^{*}_{D-1}v +FE^{*}_{D}v = \theta_{t}E_{D}^{*}v. }$
\end{item}

\end{enumerate}
\end{lemma}

\begin{proof}
    Observe that $A v = \theta_t v$.  Fix an integer $i$ 
    $(\intv0iD)$.  Now by 
    (\ref{sumrfl}), (\ref{hick}), we have
\begin{eqnarray}
    R E_{i-1}^* v  + FE_{i}^{*}v + L E_{i+1}^* v &=& E^*_i A v  \nonumber \\
                                    &=& \theta_t E^*_i v. \nonumber
\end{eqnarray}
     Assertion (i)
     follows since $FE^*_i=0$ for $0 \leq i \leq D-1$.
     Assertion (ii) similarly follows since $E^*_{D+1}=0$.
\end{proof}

\begin{lemma} \label{2ndrel}
     With reference to Definition \ref{setup},
    let $W$ denote an irreducible $T$-module with dual endpoint $t$, and fix
    any nonzero $v \in E_{t} W$.  Suppose $v$ is an eigenvector for $F^*$
    with eigenvalue $\alpha$.  Then
\begin{enumerate}
\begin{item}
$\displaystyle{\theta^*_{i-1} R E^*_{i-1} v + \theta^*_{i+1} L E^*_{i+1} v =
(\theta_{t+1} \theta^*_i - \alpha \theta_{t+1} + \alpha \theta_{t}) E^*_i v
    \hspace{1em} (\interval0i{D-1}),}$ 
\end{item}

\begin{item}
    $\displaystyle{ \theta_{D-1}^{*}RE_{D-1}^{*}v + \theta_{D}^{*}FE_{D}^{*}v = 
     (\theta_{t+1} \theta_{D}^{*} - \alpha \theta_{t+1} + \alpha 
     \theta_{t})E_{D}^{*}v,}$
\end{item}
\end{enumerate}
     where $\theta_{D+1}$, $\theta^*_{-1}$ are
    indeterminates.
\end{lemma}
 
\begin{proof}
    Observe $L^*v =0$ by (\ref{hick-s}), and
    $F^*v = \alpha v$ by assumption, so $R^*v = (A^* - \alpha I)v$
    in view of (\ref{sumrfl*}).
    Since $R^*v \in E_{t+1}W$ by (\ref{hick-s}),
\begin{equation}\label{his2}
    A (A^* - \alpha I)v = \theta_{t+1} (A^* - \alpha I) v.
\end{equation}
    Fix an integer $i$ $(\intv0iD)$.  We may now argue that
\begin{eqnarray}
\theta^*_{i-1} R E^*_{i-1}v + \theta_{i}^{*}FE_{i}^{*}v + \theta^*_{i+1} 
     L E^*_{i+1} v
    &=& (R E^*_{i-1} + FE_{i}^{*}v + L E^*_{i+1}) A^* v
    \hspace{4.9em} \mbox{(by (\ref{eieigsp*}))} \nonumber \\
    &=& E^*_i A A^* v
    \hspace{10em} \mbox{(by (\ref{sumrfl}),(\ref{hick}))} \nonumber \\
    &=& E^*_i A (A^* - \alpha I) v + \alpha E^*_i A v
    \nonumber \\
    &=& \theta_{t+1} E^*_i (A^* - \alpha I) v + \theta_{t} \alpha E^*_i v
    \hspace{2em} \mbox{(by (\ref{his2}))} \nonumber \\
    &=& (\theta_{t+1} (\theta^*_i - \alpha ) + \alpha \theta_{t}) E^*_i v
    \hspace{3.4em} \mbox{(by (\ref{eieigsp*})),} \nonumber
\end{eqnarray}
    where $\theta_{D+1}$, $\theta^*_{-1}, \theta_{D+1}^{*}$ are 
    indeterminates.  Assertion (i)
    follows since $FE^*_i=0$ for $0 \leq i \leq D-1$.
    Assertion (ii) similarly follows since $E^*_{D+1}=0$.
\end{proof}

%
%
%
%

\begin{lemma}\label{his3.8new}
    With reference to Definition \ref{setup},
    let $W$ denote an irreducible $T$-module.
    Then
    $W$ is thin and dual thin.
\end{lemma}

\begin{proof}  
    Let $t$ denote the dual endpoint of $W$.  Since $F^* E_t W
    \subseteq E_t W$, the space $E_t W$ contains a nonzero eigenvector $v$
    for $F^*$.  By Lemma \ref{1strel},
\begin{eqnarray}
    R E^*_{i-1} v + L E^*_{i+1} v & \in & span(E^*_i v)
    \hspace{3em} (\intv0i{D-1}), \label{frogface1} \\
    RE_{D-1}^{*}v + FE_{D}^{*}v & \in & span(E^{*}_{D}v).
    \label{frogface1D}
\end{eqnarray}
    By Lemma \ref{2ndrel},
\begin{eqnarray}
    \theta^*_{i-1} R E^*_{i-1} v + \theta^*_{i+1} L E^*_{i+1} v
    & \in & span(E^*_i v) \hspace{3em} (\intv0i{D-1}), 
    \label{frogface2} \\
    \theta_{D-1}^{*}RE_{D-1}^{*}v + \theta_{D}^{*}FE_{D}^{*}v & \in &
    span(E_{D}^{*}v), \label{frogface2D}
\end{eqnarray}
    where $\theta^*_{-1}$ is indeterminate.  By
    (\ref{frogface1}), (\ref{frogface2}), and (\ref{thsdist}),  we 
    find
\begin{displaymath}
    R E^*_{i} v  \in  span(E^*_{i+1} v)
    \hspace{3em} (\intv0i{D-2}), \qquad \qquad
    L E^*_{i} v  \in  span(E^*_{i-1} v) \hspace{3em} (\intv1iD).
\end{displaymath}
    By (\ref{frogface1D}), (\ref{frogface2D}), and (\ref{thsdist}), 
    we find $RE_{D-1}^{*}v \in span(E_{D}^{*}v)$ and $FE_{D}^{*}v \in 
    span(E_{D}^{*}v)$.  Combining the above information with Lemma 
    \ref{flat}(iii) and recalling 
    that $RE_{D}^{*}=0$, $LE_{0}^{*}=0$, we 
    see
\begin{eqnarray}
     R E^*_{i} v  &\in &  span(E^*_{i+1} v)
    \hspace{3em} (\intv0i{D}), \label{tadr} \\
    F E^{*}_{i}v & \in & span(E^{*}_{i}v) \hspace{3em} (\intv0iD), 
    \label{tadf} \\
    L E^*_{i} v  & \in&  span(E^*_{i-1} v) \hspace{3em} (\intv0iD).
        \label{tadl}
\end{eqnarray}
    We claim that
\begin{equation}\label{frogface5}
    W = span \{ E^*_0 v, E^*_1 v, \ldots , E^*_D v \}.
\end{equation}
    To see this, let $W'$ denote the right side of (\ref{frogface5}).
    Certainly $W' \subseteq W$; to prove that $W' = W$, we show that
    $W'$ is a nonzero $T$-module.  Observe that
    $v = \sum_{i=0}^D E^*_i v \in W',$ so $W' \not= 0$.
    Observe that $M^* W' \subseteq W'$ by
    the construction.  Observe that
    $R W' \subseteq W'$, $F W' \subseteq W'$, and $L W' \subseteq W'$ by
    (\ref{tadr})--(\ref{tadl}).  Recall that
    $A = R + F + L$ generates $M$, so $M W' \subseteq W'$.
    Since $M$, $M^*$ generate $T$, we now have that $T W' \subseteq W'$,
    so $W'$ is a $T$-module. It follows that
    $W' = W$ by the irreducibility of $W$.  We now have
    (\ref{frogface5}), which implies
    that $W$ is thin. By Lemma \ref{ppolythin},
    $W$ is dual thin.
\end{proof}


\bigskip

    We conclude this section with a comment.

\begin{theorem}\label{alphalem}
    With reference to Definition \ref{setup},
    let $W$ denote an irreducible $T$-module with diameter $d \geq 1$, 
    endpoint $r$,
    and dual endpoint $t$.  Then 
    \begin{equation} \label{alpeq}
    a^*_0(W) = \frac{\theta^*_{r+1} \theta_t -
            \theta_{t+1} \theta^*_r}{\theta_t - \theta_{t+1}}.
    \end{equation}   
\end{theorem}

\begin{proof}
     Fix any nonzero $v \in E_t W$.  Setting $i=r$ in the equation in 
     Lemma \ref{1strel}(i), we find that
\begin{equation}\label{klap1}
    L E^*_{r+1} v = \theta_t E^*_r v.
\end{equation}
    By Lemma \ref{his3.8new} and (\ref{aswiparteq}),
    $v$ is an eigenvector for $F^*$, with eigenvalue
    $a^*_0(W)$.  Setting $i=r$, $\alpha = a^*_0(W)$ in the 
    equation in Lemma \ref{2ndrel}(i), we find
\begin{equation}\label{klap2}
    \theta^*_{r+1} L E^*_{r+1} v = (\theta_{t+1} \theta^*_r -
    a^*_0(W)(\theta_{t+1}-\theta_t)) E^*_r v.
\end{equation}
    Eliminating $LE^*_{r+1}v$ in (\ref{klap2}) using (\ref{klap1}),
    and since $E^*_r v \not= 0$, we obtain
$$
    \theta^*_{r+1} \theta_t - \theta_{t+1} \theta^*_r = a^*_0(W)
        (\theta_{t}-\theta_{t+1}),
$$
    and (\ref{alpeq}) follows.
\end{proof}

%
%
%
%

\newsection{Computation of $c_i(W), \; a_{i}(W), \; b_i(W)$}

    Let $Y$ be as in Definition \ref{setup}, and
    let $W$ denote an irreducible $T$-module with diameter $d$.  In
    this section, we compute the parameters $c_i(W)$, $a_{i}(W)$, $b_i(W)$
    $(\intv0id)$.  We begin with $a_{i}(W)$.  
    
\begin{lemma}  \label{aiw0}
    With reference to Definition \ref{setup},
    let $W$ denote an irreducible $T$-module with
    diameter $d$.  Then
 \begin{enumerate}
 \begin{item}
     $a_{i}(W) =0$ $\;\;(\intv0i{d-1})$,
 \end{item}
 \begin{item}
     $a_{d}(W) \not= 0$.
 \end{item}
 \end{enumerate}
 \end{lemma}
 
\begin{proof}  (i).  Immediate from (\ref{awiparteq}) and Lemma 
\ref{flat}(iii).

\noindent (ii).  Immediate from Theorem \ref{width} and 
\cite[Theorem 15.2]{collins}.
\end{proof}

\begin{lemma}\label{cwibwieqs}
    With reference to Definition \ref{setup},
    let $W$ denote an irreducible $T$-module with
    dual endpoint $t$ and diameter $d$.
    Then
\begin{enumerate}
\begin{item}
    $\displaystyle{c_i(W) + b_i(W) = \theta_t \hspace{1em} 
    (\interval0i{d-1}),}$
\end{item}
\begin{item}
     $\displaystyle{  c_{d}(W) + a_{d}(W) = \theta_{t}.}$
\end{item}
\end{enumerate}
\end{lemma}
\begin{proof}
    Fix any nonzero $v \in E_{t} W$, and fix an integer $i$ 
    $(\intv0id)$. Using 
    (\ref{cwiparteq})--(\ref{bwiparteq}),
    Lemma \ref{1strel}, and Lemma \ref{flat}(iii), we find
\bea
    (c_i(W) + a_{i}(W) + b_i(W)) E^*_{r+i}v &=&
            (R E^*_{r+i-1} + FE_{r+i}^{*} + L E^*_{r+i+1})v \nonumber \\
        &=&\theta_t E^*_{r+i} v. \nonumber
\eea
    Since
    $E^*_{r+i} v \not= 0$ by Theorem \ref{basis}(i), we have
\begin{displaymath}
     c_{i}(W) + a_{i}(W) + b_{i}(W) = \theta_{t} \hspace{3em} (\intv0id).
\end{displaymath}
     We observe $a_{i}(W) =0 \quad (\intv0i{d-1})$ by Lemma \ref{aiw0}
     and $b_{d}(W) =0$ 
     by Definition \ref{intfunny8}.  The result follows.  
\end{proof}

\medskip

\begin{lemma}\label{cwibwieqs2}
    With reference to Definition \ref{setup},
    let $W$ denote an irreducible $T$-module with endpoint $r$,
    dual endpoint $t$, and 
    diameter $d$.
    Suppose $d \geq 1$. Then
\begin{enumerate}
\begin{item}
    $\displaystyle{\theta^*_{r+i-1} c_i(W) + \theta^*_{r+i+1} b_i(W)
        = \theta_t \theta^*_{r+1} + \theta_{t+1} \theta^*_{r+i} -
          \theta_{t+1} \theta^*_r \hspace{1em} (\intv0i{d-1}),}$
\end{item}
\begin{item}
    $\displaystyle{\theta^{*}_{r+d-1}c_{d}(W) + \theta^{*}_{r+d}a_{d}(W) = 
    \theta_t \theta^*_{r+1} + \theta_{t+1} \theta^*_{r+d} -
          \theta_{t+1} \theta^*_r,}$
\end{item}
\end{enumerate}
    where $\theta^*_{-1}$ is indeterminate.
\end{lemma}

\begin{proof}
  Fix any nonzero $v \in E_{t} W$, and fix any integer $i$ 
  $(\intv0id)$.  
    By (\ref{aswiparteq}),  $v$ is an eigenvector for $F^*$ with
    eigenvalue $a^*_0(W)$.  Using 
    (\ref{cwiparteq})--(\ref{bwiparteq}), Lemma \ref{flat}(iii),
    Lemma \ref{2ndrel}, 
    and (\ref{alpeq}), we find
\bea 
    (\theta^*_{r+i-1} c_i(W) + \theta^*_{r+i} a_i(W) +
     \theta^*_{r+i+1} b_i(W)) E^*_{r+i}v
    &=& (\theta^*_{r+i-1} R E^*_{r+i-1} + \theta^{*}_{r+i}FE_{r+i}^{*} +
            \theta^*_{r+i+1} L E^*_{r+i+1})v \nonumber \\
    &=& ( \theta_{t+1} \theta^*_{r+i} - a^*_0(W) \theta_{t+1} +
    a^*_0(W) \theta_t ) E^*_{r+i} v \nonumber \\
    &=& ( \theta_t \theta^*_{r+1} + \theta_{t +1} \theta^*_{r+i} -
    \theta_{t+1} \theta^*_{r} ) E^*_{r+i} v, \nonumber
\eea
    where $\theta_{-1}^{*}, \theta_{D+1}^{*}$ are indeterminates.
    Recall $E^*_{r+i}v \not= 0$ by
    Theorem \ref{basis}(i).  We observe $a_{i}(W) =0 \quad (\intv0i{d-1})$ 
    by Lemma 
     \ref{aiw0}(i)
     and $b_{d}(W) =0$ 
     by Definition \ref{intfunny8}.  The result follows.   
\end{proof}

\medskip


\begin{theorem}\label{cwbw-theta}
    With reference to Definition \ref{setup},
    let $W$ denote an irreducible $T$-module with endpoint $r$,
    dual endpoint $t$, and 
    diameter $d$.  First assume $d=0$.  Then
    $c_{0}(W) =0$, $a_{0}(W) = \theta_{t}$, and $b_{0}(W)=0$.
    Now assume $d \geq 1$.  Then
    \begin{eqnarray}
    c_0 (W) &=& 0, \label{c0w} \\
    c_i (W) &=& \frac{ \theta_t ( \theta^*_{r+i+1} - \theta^*_{r+1})
                     - \theta_{t+1} ( \theta^*_{r+i} - \theta^*_r )}
                   {\theta^*_{r+i+1} - \theta^*_{r+i-1}} \label{cwieq1}
                                \hspace{2em} (\interval1i{d-1}), \\
    c_d (W) &=& \frac{\theta_{t}(\theta^{*}_{r+d}-\theta^{*}_{r+1}) - 
    \theta_{t+1}(\theta^{*}_{r+d}-\theta_{r}^{*})}
    {\theta^{*}_{r+d}-\theta_{r+d-1}^{*}}, \label{cdw} \\
    \nonumber \\
    a_{i}(W) &=& 0 \hspace{2em} (\intv0i{d-1}), \label{aiw} \\
    a_{d}(W) &=& \frac{\theta_{t}(\theta^{*}_{r+d-1}-\theta^{*}_{r+1}) - 
    \theta_{t+1}(\theta^{*}_{r+d}-\theta_{r}^{*})}
    {\theta^{*}_{r+d-1}-\theta_{r+d}^{*}}, \label{adw} \\
    \nonumber \\
    b_0 (W) &=& \theta_t, \label{b0w} \\
    b_i (W) &=& \frac{ \theta_t ( \theta^*_{r+i-1} - \theta^*_{r+1})
                     - \theta_{t+1} ( \theta^*_{r+i} - \theta^*_r )}
                   {\theta^*_{r+i-1} - \theta^*_{r+i+1}} \label{bwieq1}
                                \hspace{2em} (\interval1i{d-1}), \\
    b_d (W) &=& 0. \label{bdw}
    \end{eqnarray}
    In particular, $c_i(W)$, $a_{i}(W)$, $b_i(W)$ are real for $\intv0id$.
\end{theorem}

\begin{proof}
    First assume $d=0$.  We find $a_{0}(W)=\theta_{t}$ by setting 
    $d=0$ in the equation in Corollary \ref{tracelem2}(i).  By Definition 
    \ref{intfunny8}, we find $c_{0}(W)=0$, $b_{0}(W)=0$.  
    
    Now assume 
    $d \geq 1$.  
    Lines (\ref{c0w}), (\ref{bdw}) are immediate from Definition
    \ref{intfunny8}, and line (\ref{aiw}) follows from Lemma 
    \ref{aiw0}.  Line (\ref{b0w}) follows from Lemma \ref 
    {cwibwieqs}(i) and 
    (\ref{c0w}).  
    To obtain (\ref{cwieq1}) and (\ref{bwieq1}),
    solve the linear system determined by the equations in Lemma 
    \ref 
    {cwibwieqs}(i), Lemma \ref {cwibwieqs2}(i) for the
    variables $c_i(W)$, $b_i(W)$. We observe that the coefficient
    matrix of this system is nonsingular since $\theta^*_0,$ ...,
    $\theta^*_D$ are distinct.  To obtain (\ref{cdw}) and 
    (\ref{adw}), solve the linear system determined by the equations in
    Lemma \ref{cwibwieqs}(ii), Lemma \ref{cwibwieqs2}(ii)
    for the variables $c_{d}(W)$, $a_{d}(W)$.  We observe that the coefficient
    matrix of this system is nonsingular since $\theta^*_0,$ ...,
    $\theta^*_D$ are distinct.
\end{proof}


\newsection{Computation of $c^*_i(W), \; a^*_i(W), \; b^*_i(W)$}

    Let $Y$ be as in Definition \ref{setup}, and
    let $W$ denote an irreducible $T$-module with diameter $d$.  In
    this section, we compute the parameters $c^*_i(W)$,
    $a^*_i(W)$, $b^*_i(W)$ $(\intv0id)$.

\begin{lemma}\label{cswibswieqs}
    With reference to Definition \ref{setup},
    let $W$ denote an irreducible $T$-module with endpoint $r$
    and diameter $d$.  Fix any integer $i \; (\interval0id)$. Then
    \begin{equation}\label{strl1eq}
    c^*_i(W) + a^*_i(W) + b^*_i(W) = \theta^*_r.
    \end{equation}
\end{lemma}

\begin{proof}
    Let $t$ denote the dual endpoint of $W$.
    Pick any nonzero $v \in E^*_r W$, and observe that $A^* v = \theta^*_r v$.
    We may now argue that
\bea
    (c_i^*(W)&+&a_i^*(W) \; \;  +  \; \; b_i^*(W)) E_{t+i}v  \nonumber \\
    && \; \; \; = \; (R^* E_{t+i-1} \; + \; F^* E_{t+i} \; + \; L^* E_{t+i+1})v
         \hspace{5em} (\mbox{by Def. \ref{int*}}) \nonumber \\
    && \; \; \; = \; E_{t+i} A^* v \hspace{17.2em} (\mbox{by (\ref{sumrfl*}),
    (\ref{hick-s})})
                                                    \nonumber\\
    && \; \; \; = \; \theta^*_r E_{t+i} v.  \nonumber
\eea
    The result now follows, since
    $E_{t+i} v \not= 0$ by Theorem \ref{basis}(ii).
\end{proof}

\begin{lemma}\label{cswibswieqs2}
    With reference to Definition \ref{setup},
    let $W$ denote an irreducible $T$-module with endpoint $r$,
    dual endpoint $t$, and diameter $d$.
    Suppose $d \geq 1$, and fix any integer
    $i$ $(\interval0id)$. Then
    \begin{equation}\label{strl2eq}
    \theta_{t+i-1} c^*_i(W) + \theta_{t+i} a^*_i(W)
            + \theta_{t+i+1} b^*_i (W)
        = \theta^*_{r+1} \theta_{t+i},
    \end{equation}
    where $\theta_{-1}$, $\theta_{D+1}$ are indeterminates.
\end{lemma}

\begin{proof}
    Pick any nonzero $v \in E^*_r W$. Observe $r<D$ since $d \geq 1$, 
    so $Fv=0$ by Lemma \ref{flat}(iii).  Hence 
    $Av = R v \in E^*_{r+1} W.$ It follows that
\begin{equation}\label{jdkk}
    A^* Av = \theta^*_{r+1} Av.
\end{equation}
    We may now argue that
\bea
    &&(\theta_{t+i-1} c_i^*(W) \; + \; \theta_{t+i}
            a_i^*(W) \; + \; \theta_{t+i+1}
                                        b_i^*(W)) E_{t+i}v  \nonumber \\
    &&\hspace{3em} = \; \; \; (\theta_{t+i-1} R^* E_{t+i-1} \; + \;
    \theta_{t+i}
            F^* E_{t+i} \; + \; \theta_{t+i+1} L^* E_{t+i+1})v
         \hspace{3.3em} (\mbox{by Def. \ref{int*}}) \nonumber \\
    &&\hspace{3em} = \; \; \; (R^* E_{t+i-1} +  F^* E_{t+i} + L^* E_{t+i+1}) A v
                        \hspace{10.1em} (\mbox{by (\ref{eieigsp})}) \nonumber \\
    &&\hspace{3em} = \; \; \; E_{t+i} A^* A v
                        \hspace{21.7em} (\mbox{by
        (\ref{sumrfl*}), (\ref{hick-s})}) \nonumber \\
    &&\hspace{3em} = \; \; \; \theta^*_{r+1} E_{t+i} A v
                        \hspace{21.1em} (\mbox{by (\ref{jdkk})}) \nonumber \\
    &&\hspace{3em} = \; \; \; \theta^*_{r+1} \theta_{t+i} E_{t+i} v.
                        \hspace{20.1em} (\mbox{by (\ref{eieigsp})}) \nonumber
\eea
    The result now follows, since
    $E_{t+i} v \not= 0$ by Theorem \ref{basis}(ii).
\end{proof}

\begin{lemma}\label{cswibswieqs3}
    With reference to Definition \ref{setup},
    let $W$ denote an irreducible $T$-module with endpoint $r$,
    dual endpoint $t$, and diameter $d$.
    Suppose $d \geq 2$, and fix any integer
    $i$ $(\interval0id)$. Then
    \begin{eqnarray}
    \theta^2_{t+i-1} \; c^*_i(W) \; + \; \theta^2_{t+i} \; a^*_i(W) &+&
    \theta^2_{t+i+1} \; b^*_i(W)  \nonumber \\
     = \; \theta^*_{r+2}  \theta^2_{t+i} &+&
    b_0(W)c_1(W) (\theta^*_r - \theta^*_{r+2}),\label{strl3eq}
    \end{eqnarray}
    where $\theta_{-1}$, $\theta_{D+1}$ are indeterminates.
\end{lemma}

\begin{proof}
    For notational convenience, set $\alpha := b_0(W) c_1(W)$.
    Pick any nonzero $v \in E^*_r W$.
    We first claim that $LRv = \alpha v$.
    To see this, observe by Theorem \ref{basis}(i),
    there exists a nonzero $z \in E_t W$
    such that $v = E^*_r z$. Applying Definition
    \ref{intfunny8},
\bea
    LRv &=& LR E^*_r z \nonumber \\
        &=& c_1(W) L E^*_{r+1} z \nonumber \\
        &=& b_0(W) c_1(W) E^*_r z, \nonumber
\eea
    and the claim follows.
    
    Observe $r+1 < D$ since $d \geq 2$, so $Fv=0$, $FRv=0$ by Lemma 
    \ref{flat}(iii).  By these remarks, the above claim, line (\ref{sumrfl}), 
    and since $Lv =0$,
\begin{eqnarray}
    R^2 v &=& (R+F+L)^{2} v - LRv  \nonumber \\
          &=& (A ^2 - \alpha I) v.  \label{hick2}
\end{eqnarray}
    By (\ref{hick2}), and since $R^2 v \in E^*_{r+2} V$ by
    (\ref{hick}),
\begin{equation}\label{hick3}
    A^* (A^2-\alpha I)v = \theta^*_{r+2} (A^2 - \alpha I)v .
\end{equation}
    We may now argue that
\bea
    &&(\theta^2_{t+i-1} c_i^*(W) \; + \; \theta^2_{t+i}
            a_i^*(W) \; + \; \theta^2_{t+i+1}
                                        b_i^*(W)) E_{t+i}v  \nonumber \\
    &&\hspace{3em}= \; \; \; (\theta^2_{t+i-1} R^* E_{t+i-1} \; +
        \; \theta^2_{t+i}
            F^* E_{t+i} \; + \; \theta^2_{t+i+1} L^* E_{t+i+1})v
         \hspace{3.3em} (\mbox{by Def. \ref{int*}}) \; \; \; \nonumber \\
    &&\hspace{3em}= \; \; \; (R^* E_{t+i-1} +  F^* E_{t+i} + L^* E_{t+i+1})
    A^2 v  \hspace{10em} (\mbox{by (\ref{eieigsp})}) \nonumber \\
    &&\hspace{3em}= \; \; \; E_{t+i} A^* A^2 v
                        \hspace{21em} (\mbox{by
        (\ref{sumrfl*}),(\ref{hick-s})}) \nonumber \\
    &&\hspace{3em}= \; \; \; E_{t+i} A^* (A^2 - \alpha I) v  +
        \alpha E_{t+i} A^* v
                        \nonumber \\
    &&\hspace{3em}= \; \; \; \theta^*_{r+2} (A^2 - \alpha I) E_{t+i} v +
         \alpha \theta^*_r E_{t+i} v
                        \hspace{12.1em} (\mbox{by (\ref{hick3})}) \nonumber \\
    &&\hspace{3em}= \; \; \; (\theta^*_{r+2} (\theta_{t+i}^2 - \alpha) +
         \alpha \theta^*_r) E_{t+i} v.
                        \hspace{13em} (\mbox{by (\ref{eieigsp})}). \nonumber
\eea
    The result now follows, since
    $E_{t+i} v \not= 0$ by Theorem \ref{basis}(ii).
\end{proof}


\begin{theorem}\label{cswbsw-theta}
    With reference to Definition \ref{setup},
    let $W$ denote an irreducible $T$-module with endpoint $r$,
    dual endpoint $t$, and diameter $d$.  First assume $d=0$.  Then
    $c_{0}^{*}(W)=0$, $a_{0}^{*}(W) = \theta_{r}^{*}$, $b_{0}^{*}(W) 
    =0$.  Now assume $d \geq 1$.  Then
\begin{eqnarray}
    c^*_0 (W) &=& 0, \label{c0w*} \\
    c^*_i (W) &=& \frac{ ( \theta^2_{t+i} - \theta^2_t)
                        ( \theta^*_{r+2} - \theta^*_{r+1} ) +
                (\theta_t \theta_{t+1} - \theta_{t+i} \theta_{t+i+1} )
                        ( \theta^*_{r+1} - \theta^*_r ) }
    {(\theta_{t+i-1} - \theta_{t+i})(\theta_{t+i-1} -
    \theta_{t+i+1})} \hs3
            (\interval1i{d-1}), \label{cswieq1}  \\
    c^*_d (W) &=& \frac{\theta_{t+d} (\theta^*_{r+1}-\theta^*_{r})}
                             {\theta_{t+d-1} - \theta_{t+d}}, \label{cdw*} \\
    \nonumber \\
    b^*_0 (W) &=& \frac{\theta_t (\theta^*_r-\theta^*_{r+1})}
                             {\theta_t - \theta_{t+1}}, \label{b0w*} \\
    b^*_i (W) &=& \frac{ ( \theta^2_{t+i} - \theta^2_t)
                        ( \theta^*_{r+2} - \theta^*_{r+1} ) +
                (\theta_t \theta_{t+1} - \theta_{t+i} \theta_{t+i-1} )
                        ( \theta^*_{r+1} - \theta^*_r ) }
    {(\theta_{t+i+1} - \theta_{t+i})(\theta_{t+i+1} -
    \theta_{t+i-1})} \hs3
            (\interval1i{d-1}), \label{bswieq1} \\
    b^*_d (W) &=& 0. \label{bdw*} \\
    \nonumber \\
    a^*_i (W) &=& \theta^*_r - b^*_i(W) - c^*_i(W) \hspace{3em}
                            (\interval0id). \label{aiw*}
    \end{eqnarray}
    In particular, $c^*_i(W)$, $a^*_i(W)$, $b^*_i(W)$ are real for $\intv0id$.
\end{theorem}

\begin{proof}
    First assume $d=0$.  We find $a_{0}^{*}(W)=\theta_{r}^{*}$ by setting 
    $d=0$ in the equation in Corollary \ref{tracelem2}(ii).  By Definition 
    \ref{int*}, we find $c_{0}^{*}(W)=0$, $b_{0}^{*}(W)=0$.  
    
    Now assume 
    $d \geq 1$.
    Lines (\ref{c0w*}), (\ref{bdw*}) follow from Definition 
    \ref{int*}.  
    Line (\ref{aiw*}) follows from (\ref{strl1eq}).
    We obtain (\ref{cdw*}) by
    setting $i=d$ in (\ref{strl1eq}), (\ref{strl2eq}), and solving
    for $c^*_d(W)$, using (\ref{bdw*}). We now have (\ref{cdw*}), and
    line (\ref{b0w*}) is obtained similarly.

    It remains to prove (\ref{cswieq1}) and (\ref{bswieq1}).
    Assume $d \geq 2$ and fix any $i \; (\interval1i{d-1}).$
    Observe (\ref{strl1eq}), (\ref{strl2eq}), and (\ref{strl3eq})
    form a system of three linear equations in the three variables
    $c^*_i(W)$, $a^*_i(W)$, and $b^*_i(W)$.
    Observe the coefficient matrix is Vandermonde,
    hence nonsingular, since $\theta_0$,...,$\theta_D$
    are distinct.  Solving this system, we
    obtain (\ref{cswieq1}) and (\ref{bswieq1}).
\end{proof}


\newsection{Some Comments on the Intersection Numbers}

    Let $Y$ be as in Definition \ref{setup}, and
    let $W$ denote an irreducible $T$-module with diameter $d$.  
    In this section,
    we show the expressions $b_{i-1}(W)c_{i}(W)$ and 
    $ b_{i-1}^{*}(W)c_{i}^{*}(W)$ are 
    positive for $\intv1id$.  

\begin{lemma}\label{uiidd}
    With reference to Definition \ref{setup},
    let $W$ denote
    an irreducible $T$-module with endpoint $r$, dual endpoint $t$,
    and diameter $d$.
\begin{enumerate}
\begin{item} \label{sqrnrmsI}
    For any nonzero $v \in E_{t} W$,
\begin{displaymath}
    c_i(W) \| E^*_{r+i} v \| ^2 = b_{i-1}(W) \| E^*_{r+i-1} v \| ^2
    \hspace{3em} (\intv1id).
\end{displaymath}
\end{item}
\begin{item} \label{sqrnrmI}
    For any nonzero $v \in E^*_{r} W$,
\begin{displaymath}
    c^*_i(W) \| E_{t+i} v \| ^2 = b^*_{i-1}(W) \| E_{t+i-1} v \| ^2
    \hspace{3em} (\intv1id).
\end{displaymath}
\end{item}
\end{enumerate}
\end{lemma}
\begin{proof}
(i).  By (\ref{cwiparteq}), (\ref{bwiparteq}),
    and since $R = \overline{L}^t$,
\begin{eqnarray}
      c_i(W) \; \| E^*_{r+i} v \| ^2 &=&
    \langle R E^*_{r+i-1}v , \; E^*_{r+i}v \rangle  \nonumber \\
    &=& \langle E^*_{r+i-1}v , \; L E^*_{r+i}v \rangle  \nonumber \\
                    &=& \overline{b_{i-1}(W)} \; \| E^*_{r+i-1} v \| ^2.
    \nonumber
\end{eqnarray}
    Recall that $b_{i-1}(W)$ is real by Theorem \ref{cwbw-theta}, so the
    result follows.

\noindent
(ii).  Similar to the proof of (i).
\end{proof}

\begin{corollary}\label{8394}
    With reference to Definition \ref{setup},
    let $W$ denote an irreducible $T$-module with diameter $d$.
    Then
\begin{enumerate}
\begin{item}
$b_{i-1}(W)c_i(W) > 0 \hspace{3em} (\intv1id)$.
\end{item}
\begin{item}
$b^*_{i-1}(W)c^*_i(W) > 0  \hspace{3em} (\intv1id)$.
\end{item}
\end{enumerate}
\end{corollary}

\begin{proof}
(i).  Let $r$ denote the endpoint of $W$.  
    The product $b_{i-1}(W) c_{i}(W)$ is nonnegative by Lemma
    \ref{uiidd}(i), and since $\| E^*_{r+i} v \| ^2$ and
    $\| E^*_{r+i-1} v \| ^2$ are positive.  Observe
    $b_{i-1}(W) \not= 0$ by (\ref{bwiparteq}) and the 
    previously-mentioned fact that $LE_{i}^{*}W \not=0$ $(r < i  \leq 
    r+d)$.  
    Similarly
    $c_{i}(W) \not= 0$ by (\ref{cwiparteq}) and the 
    fact that $RE_{i}^{*}W \not=0$ $(r \leq i  < 
    r+d)$.  

\noindent
(ii).  Similar to the proof of (i).
\end{proof}


\newsection{The Isomorphism Classes of Irreducible $T$-Modules}

    With reference to Definition \ref{setup},
    in this section we prove that the isomorphism
    class of any given irreducible $T$-module 
    is determined by its dual endpoint and diameter.

\bigskip

\begin{theorem}\label{isom}
    With reference to Definition \ref{setup},
    let $W$, $W'$ denote irreducible $T$-modules with endpoints $r$,
    $r'$, dual endpoints $t$, $t'$, and diameters $d$, $d'$, respectively.
    Then the following are equivalent.
\begin{enumerate}
\begin{item}\label{isompart1}
    $W$ and $W'$ are isomorphic as $T$-modules.
\end{item}
\begin{item}\label{isompart2}
    $t = t'$ and $d = d'$.
\end{item}
\begin{item}\label{isompart3}
    $t = t'$ and $r = r'$.
\end{item}
\begin{item}\label{isompart4}
    $B(W) = B(W')$.
\end{item}
\begin{item}\label{isompart5}
    $t = t'$ and $B^{*}(W) = B^{*}(W')$.
\end{item}
\end{enumerate}
\end{theorem}

\begin{proof}
(\ref{isompart1}) $\rightarrow$ (\ref{isompart2}).
    Let $\phi$ denote a $T$-isomorphism from $W$ to $W'$.
    Then for any integer $i \; (\interval0iD)$,
\begin{eqnarray}
    &E_i W =0 \; \Leftrightarrow \; \phi (E_i W) =0 \; \Leftrightarrow \;
         E_i \phi (W) =0 \; \Leftrightarrow \; E_i W' =0. \nonumber
\end{eqnarray}
    Now (\ref{isompart2}) follows by (\ref{d*defeq}) and since $d=d^{*}$.

\noindent
(\ref{isompart2}) $\leftrightarrow$ (\ref{isompart3}).
    Immediate by Theorem \ref{width}.

\noindent
(\ref{isompart2}),(\ref{isompart3}) $\rightarrow$ (\ref{isompart4}).
    By Theorem \ref{cwbw-theta}, the entries in the intersection
    matrix are determined by the endpoint, dual endpoint, and diameter.
    
\noindent
(\ref{isompart2}),(\ref{isompart3}) $\rightarrow$ (\ref{isompart5}).
    By Theorem \ref{cswbsw-theta}, the entries in the dual intersection
    matrix are determined by the endpoint, dual endpoint, and diameter.

\noindent
(\ref{isompart4}) $\rightarrow$ (\ref{isompart2}).
    $B(W)$ is a $d+1 $ by $\;d+1 $ matrix, and
    $B(W')$ is $\; d'+1 $ by $\;d'+1, $ so $d=d'$.
    By Lemma \ref{aiw0} and Lemma \ref{cwibwieqs}, the sum  
    of the entries in each
    row of $B(W)$ equals $\theta_{t}$, and the sum of the entries in 
    each row of $B(W')$ equals $\theta_{t'}$.
    Hence $\theta_{t} = \theta_{t'}$, so $t=t'$ in
    view of (\ref{thdist}).

\noindent
(\ref{isompart5}) $\rightarrow$ (\ref{isompart1}).
    $B^*(W)$ is a $\; d+1 $ by $ d+1 \;$ matrix, and
    $B^*(W')$ is $\; d'+1 $ by $ d'+1 \;$, so $d=d'$.
    Now $r = r'$ by Theorem \ref{width}.  Pick a
    nonzero  $v \in E_r^{*} W$ and recall by
    Theorem \ref{basis}(ii) that
    $ E_t v, \ldots , E_{t+d}v $ is a basis
    for $W$.
    Similarly, pick a nonzero $v' \in E_r^{*} W'$, and observe
    $E_t v', \ldots , E_{t+d}v' $ is a basis for $W'$.
    By linear algebra, there exists an isomorphism of vector spaces
    $\phi : W \rightarrow W'$ such that
\begin{equation} \label{phidef}
    \phi \, : \; E_i v \mapsto E_i v' \hspace{3em} (\interval{t}i{t+d}).
\end{equation}
    We show $\phi$ is an isomorphism of $T$-modules.
    Since $A^{*}, E_0, \ldots, E_D$ generate $T$,
    and since $A^{*}=R^{*}+F^{*}+L^{*}$, it suffices to show
\begin{eqnarray}
    (R^{*} \phi - \phi R^{*})W &=& 0, \label{Req} \\
    (F^{*} \phi - \phi F^{*})W &=& 0, \label{Feq} \\
    (L^{*} \phi - \phi L^{*})W &=& 0, \label{Leq} \\
    (E_j \phi - \phi E_j)W &=& 0 \hspace{2.5em}
        (\interval0{j}D). \label{Ejeq}
\end{eqnarray}
    Line (\ref{Ejeq}) is immediate from the construction.
    To see (\ref{Req}), observe that $c^{*}_i(W) = c^{*}_i(W')$ 
    $(\interval0id)$.
    Now by (\ref{cswiparteq}), $R^{*} \phi - \phi R^{*}$ vanishes on 
    $E_{t+i}v$ $(\intv0id)$.  Line (\ref{Req}) follows.   
    Lines (\ref{Feq}), (\ref{Leq}) are proven similarly.
\end{proof}


\bigskip

    We conclude this section with a few comments.

\begin{lemma}\label{trivmod}  \cite[Theorem 4.1]{terwSubconstituentII}
    With reference to Definition \ref{setup}, there exists
a unique irreducible $T$-module $W_0$ with diameter $D$.
The endpoint and dual endpoint of $W_0$ are both zero.
We refer to $W_0$ as the {\em trivial $T$-module}.
\end{lemma}

\begin{lemma}\label{trivthm1}\cite[Theorem 4.1]{terwSubconstituentII}
    With reference to Definition \ref{setup},  let $W_0$ denote
    the trivial module for $Y$. Then
\begin{enumerate}
\begin{item} $\; \; c_i(W_0) = c_i, \hspace{2em} a_i(W_0) = a_i, \hspace{2em} 
b_i(W_0) = b_i
     \hspace{3em} (\interval0iD)$,
\end{item}
\begin{item} $\; \; c^*_i(W_0) = c^*_i, \hspace{2em} a^*_i(W_0) = a^*_i,
    \hspace{2em} b^*_i(W_0) = b^*_i \hspace{3em} (\interval0iD)$.
\end{item}
\end{enumerate}
\end{lemma}


\newsection{The Parameters in Terms of $q$ and $s$}

    With reference to Definition \ref{setup}, we 
    now explicitly compute the intersection
    matrices and dual intersection matrices
    of the irreducible $T$-modules.  For convenience, we exclude a 
    small class of examples.  
    
    Let $D$ denote any integer at least 3.  Let $O_{D+1}$ denote the 
    Odd graph with diameter $D$, and let $\Box_{2D+1}$ denote the folded 
    $(2D+1)$-cube \cite[pp. 259, 264]{bcn}.  It is well-known that $O_{D+1}$, $\Box_{2D+1}$ are
    almost-bipartite P- and Q-polynomial schemes with diameter $D$.  
    In \cite{moon}, it was shown that $O_{D+1}$ is uniquely determined by 
    its intersection numbers.  By \cite[p. 264]{bcn},  
    $\Box_{2D+1}$
    is uniquely determined by its intersection numbers.  For more 
    information on the structure of the irreducible $T$-modules for 
    these schemes, see \cite{terwSubconstituentIII}.

\begin{lemma}\label{qbipform}
    With reference to Definition \ref{setup},  suppose
    $Y$ is not one of $O_{D+1}$, $\Box_{2D+1}$.  
    Then there exist scalars $q$, $s$, $h$, $h^* \in \complexes$, 
    with $q$, $h$, $h^*$ nonzero, such that
\begin{eqnarray}
    \theta_i &=& \theta_{0}+h(1-q^{i})(1-sq^{i+1})q^{-i} \hspace{1em} (\intv0iD)
    , \label{du9i} \\
    \theta^*_i &=&  \theta^*_0 + h^*(1-q^i)(1-q^{-2D-1+i})q^{-i}
    \hspace{1em} (\interval0iD). \label{thsiqs}
\end{eqnarray}
\end{lemma}
\begin{proof}
    By \cite[pp. 237,240]{bcn}, there exist scalars
    $\beta, \gamma, \gamma^* \in \reals$ such that
\begin{eqnarray}
    \theta_{i-1} - \beta \theta_i + \theta_{i+1} &=& \gamma
                                    \hspace{3em} (\intv1i{D-1}),
                                                    \label{threcur}\\
    \theta^*_{i-1} - \beta \theta^*_i + \theta^*_{i+1} &=& \gamma^*
                                    \hspace{3em} (\intv1i{D-1}).
                                                    \label{thsrecur}
\end{eqnarray}

    We show $\beta \not= 2$, $\beta \not= -2$.  
    First, suppose $\beta=2$.  By \cite[Theorem 1.11.1]{bcn}, there exists a bipartite 
    P-polynomial scheme $Y'$ whose quotient scheme is $Y$.  
    Applying 
    \cite[Theorems 10.4, 15.2]{lang} to $Y'$, we find $Y$ is 
    $\Box_{2D+1}$, a contradiction.
    Therefore $\beta \not= 2$.

    Now suppose $\beta=-2$.  Then by \cite[Theorem 2]{terwqminus1}, $Y$ must be 
    $O_{D+1}$ or $\Box_{2D+1}$, a contradiction.  
    Therefore $\beta \not= -2$.  

    Since $\beta \not= 2$  and $\beta \not= -2$,
    there exists $q \in \complexes$ such that
    $q \not\in \{ 1, 0, -1 \}$, and
    such that $\beta = q + q^{-1}$.
    Solving the recurrence in (\ref{threcur}), we obtain
\begin{equation}\label{39ikj}
    \theta_i = h q^{-i} + h' q^i +h'' \hspace{3em} (\intv0iD)
\end{equation}
    for some $h, h', h'' \in \complexes$.
    By (\ref{thdist}), $h, h^{\prime}$ are not both zero.
    Replacing $q$ by $q^{-1}$ if necessary, we may
    assume $h \not= 0$.  Now there exists $s \in
    \complexes$ such that $h^{\prime} = s q h$.
    Eliminating $h^{\prime}$ in (\ref{39ikj}) using this, we find
\begin{displaymath}
  \theta_i - \theta_0 =   h(1-q^i)(1-sq^{i+1})q^{-i}
    \hspace{3em} (\intv0iD), 
    \end{displaymath}
    and (\ref{du9i}) follows.

    Solving (\ref{thsrecur}), we obtain
\begin{equation}\label{39ikjT}
    \theta^*_i = h^* q^{-i} + h^{*\prime} q^i + h^{*\prime\prime}
    \hspace{3em} (\intv0iD)
\end{equation}
    for some $h^*, h^{*\prime}, h^{*\prime\prime} \in
    \complexes$.  By 
    \cite[Theorem 2]{terwPQAntipodal} we have $h^{*\prime} =  h^* q s^{*}$, where
    $s^{*}=q^{-2D-2}$.  Eliminating $h^{*\prime}$ in (\ref{39ikjT}) 
    using this, we obtain (\ref{thsiqs}).  Observe $h^{*}$ is nonzero 
    by  (\ref{thsdist}).  
\end{proof}

\begin{corollary}\label{nonzeros}
    With reference to Definition \ref{setup},  suppose
    $Y$ is not one of $O_{D+1}$, $\Box_{2D+1}$.
    Let $q, s$ be as in
    Lemma \ref{qbipform}.  Then 
\bea
    q^{i} & \not = & 1 \hspace{3em} (\intv1i{2D}), \\
     sq^{i} & \not= & 1 \hspace{3em} (\intv2i{2D}).
\eea
\end{corollary}
\begin{proof}
   Evaluate (\ref{thdist}), (\ref{thsdist}) in terms of $q,s$ 
   using Lemma \ref{qbipform}.
\end{proof}

\begin{lemma}\label{qbipformII}
    With reference to Definition \ref{setup},  suppose
    $Y$ is not one of $O_{D+1}$, $\Box_{2D+1}$.
    Let $q, s, h, h^*$ be as in
    Lemma \ref{qbipform}.  Then 
\begin{eqnarray}
    \theta_{0} &=& h(1+sq), \label{th0qs} \\
    \theta^*_0 &=&  \frac{h^* (q^{2D}-1)(1+sq)}{q^{2D} (1-sq^{2})}.
                                    \label{ths0qs}
\end{eqnarray}
\end{lemma}

\begin{proof}
    We apply (\ref{strl3eq}) with $W=W_{0}$, $r=0$, $t=0$, $d=D$, 
    and $i=D$.  Evaluating 
      the result using Theorem \ref{cwbw-theta}, Theorem 
      \ref{cswbsw-theta}, and Lemma \ref{qbipform}, we obtain (\ref{th0qs}).
      
    Recall that $a_{0}^{*}=0$, so by 
    (\ref{alpeq}) (with $r=0$, $t=0$),
\begin{displaymath}
    \theta_{1}^{*}\theta_{0} = \theta_{1} \theta_{0}^{*}.
\end{displaymath}
     Solving this equation for $\theta_{0}^{*}$ using (\ref{du9i}),
     (\ref{thsiqs}), (\ref{th0qs}), we obtain (\ref{ths0qs}).  
\end{proof}

\begin{corollary}  \label{newthi}
    With reference to Definition \ref{setup},  suppose
    $Y$ is not one of $O_{D+1}$, $\Box_{2D+1}$.
    Let $q, s, h, h^*$ be as in
    Lemma \ref{qbipform}.  Then    
\be   \label{newthqs}
    \theta_{i} = hq^{-i}(1+sq^{2i+1}) \hspace{3em} (\intv0iD).
\ee
\end{corollary}

\begin{proof}
    Routine using (\ref{du9i}), (\ref{th0qs}).
\end{proof}

\begin{theorem}\label{cwbw-qs*}
    With reference to Definition \ref{setup}, suppose
    $Y$ is not one of $O_{D+1}$, $\Box_{2D+1}$. 
    Let $W$ denote an irreducible $T$-module
    with dual endpoint $t$ and diameter $d$.  First assume $d=0$.  
    Then $c_{0}(W) =0$, $a_{0}(W) = hq^{-t}(1+sq^{2t+1})$, 
    $b_{0}(W)=0$, where $q$, $s$, $h$ are as in Lemma \ref{qbipform}.  
    Now assume $d \geq 1$.  Then
    \begin{eqnarray}
    c_0 (W) &=& 0,  \label{tt532} \\
    c_i (W) &=& \frac{ h 
    (1-q^{i})(1+sq^{2+2d+2t-i})}{q^{t+i}(q^{2d-2i+1}-1)} \hs3 
         (\interval1i{d-1}),  \label{ciwqs}\\
    c_d (W) &=& \frac{h(1-q^{d})(1+sq^{2+d+2t})}{q^{t+d}(q-1)},  
         \label{tt533} \\
    a_{i}(W) &=& 0 \hs3 (\intv0i{d-1}), \\
    a_{d}(W) &=& \frac{h(q^{d+1}-1)(1+sq^{1+d+2t})}{q^{t+d}(q-1)}, 
    \label{adwqs} \\
    b_0 (W) &=& hq^{-t}(sq^{2t+1}+1), \label{tt534} \\
    b_i (W) &=& \frac{ h 
    (q^{2d+1-i}-1)(1+sq^{2t+i+1})}{q^{t+i}(q^{2d-2i+1}-1)} \hs3 
         (\interval1i{d-1}),
                    \label{biwformulathing} \\
    b_d (W) &=& 0 , \label{tt535}
    \end{eqnarray}
    where $q$, $s$, $h$ are from Lemma \ref{qbipform}.
\end{theorem}

\begin{proof}
    Routine using Theorem \ref{cwbw-theta}, Theorem 
    \ref{width}, (\ref{thsiqs}), (\ref{ths0qs}), and
    (\ref{newthqs}).
\end{proof}

\begin{theorem}\label{cwbws-qs}
    With reference to Definition \ref{setup}, suppose
    $Y$ is not one of $O_{D+1}$, $\Box_{2D+1}$. 
    Let $W$ denote an irreducible $T$-module
    with endpoint $r$, dual endpoint $t$, and diameter $d$.  Then
    \begin{eqnarray}
    c^*_0 (W) &=& 0,  \\
    c^*_i (W) &=&  \frac{ h^* 
    (1-q^{2i})(1-s^{2}q^{2+2d+4t+2i})}{q^{D+d+1}(1-sq^{2i+2t})(1-sq^{1+2i+2t})} 
    \hs2 (\interval1i{d-1}),  \label{csiwqs}\\
    c^*_d (W) &=& \frac{h^* 
    (1-q^{2d})(1+sq^{2t+2d+1})}{q^{D+d+1}(1-sq^{2t+2d})},  \\
    b^*_0 (W) &=& \frac{h^* 
    (q^{2d}-1)(1+sq^{2t+1})}{q^{D+d}(1-sq^{2+2t})}, \\
    b^*_i (W) &=& \frac{ h^* 
    (q^{2d-2i}-1)(1-s^{2}q^{2+2i+4t})}
    {q^{D+d-2i}(1-sq^{2+2i+2t})(1-sq^{1+2i+2t})} \hs2 
    (\interval1i{d-1}), \\
    b^*_d (W) &=& 0, \\
    a^{*}_{i}(W) &=& \theta_{r}^{*} - b^*_i (W) - c^*_i (W)  
    \hs2 (\interval0i{d}),
    \end{eqnarray}
    where $q$, $s$, $h^{*}¥$ are from Lemma \ref{qbipform}.
\end{theorem}

\begin{proof}
    Routine using Theorem \ref{cswbsw-theta}, Theorem 
    \ref{width}, (\ref{thsiqs}), (\ref{ths0qs}), and
    (\ref{newthqs}).
\end{proof}

\begin{corollary} \label{bounds}
  With reference to Definition \ref{setup},  suppose
    $Y$ is not one of $O_{D+1}$, $\Box_{2D+1}$.
    Let $q, s$ be as in
    Lemma \ref{qbipform}.  Then 
\begin{displaymath}
     sq^{i}  \not=  -1 \hspace{3em} (\intv1i{2D+1}). 
\end{displaymath}
\end{corollary}

\begin{proof}
   Recall the trivial module $W_{0}$ has dual endpoint $t=0$ and 
   diameter $d=D$.  By Lemma \ref{trivthm1}, the intersection numbers 
   and dual intersection numbers of $W_{0}$ in (\ref{ciwqs}), 
   (\ref{adwqs})--(\ref{biwformulathing}), (\ref{csiwqs}) are 
   nonzero.  Using this information, we obtain the desired result.
\end{proof}

\begin{corollary}  \label{hhs}
     With reference to Definition \ref{setup}, suppose
    $Y$ is not one of $O_{D+1}$, $\Box_{2D+1}$.  Then
\begin{eqnarray}  
    h &=& \frac{q-q^{2D}}{(q-1)(1+sq^{2D+1})}, \label{ryeu} \\
      h^* &=& 
      \frac{q^{2D+1}(1-sq^{2})(1-sq^{3})}{(1-q^{2})(1-s^{2}q^{2D+4})},   
      \label{ryeu2}                         
\end{eqnarray}
     where $q$, $s$, $h$, $h^{*}$ are from Lemma \ref{qbipform}.  
     We remark the denominators in (\ref{ryeu}), (\ref{ryeu2}) are 
     nonzero by Corollary \ref{nonzeros} and Corollary \ref{bounds}.
\end{corollary}

\begin{proof}
    Concerning (\ref{ryeu}), recall the trivial module $W_{0}$ has 
    endpoint $r=0$, dual endpoint $t=0$, and diameter $d=D$.  Also 
    $c_{1}(W_{0})=1$ by Lemma \ref{trivthm1} and since $c_{1}=1$.  
    Evaluating (\ref{ciwqs}) using this data we find
\begin{displaymath}
    1 = \frac{h(q-1)(1+sq^{2D+1})}{q-q^{2D}},
\end{displaymath}
    and line (\ref{ryeu}) follows.  Line (\ref{ryeu2})
    is proven similarly.
\end{proof}


\newsection{Multiplicities of the Irreducible T-Modules}

With reference to Definition \ref{setup}, in this section
we compute the multiplicities with which the irreducible
$T$-modules appear in the standard module $V$.

\begin{definition}\label{multiplicitydef}
With reference to Definition \ref{setup},
fix a decomposition of the standard module $V$ into an
orthogonal direct sum of irreducible $T$-modules.
For any integers $t,d \; (\intv0{t,d}D)$, we define $mult(t,d)$ to be
the number of irreducible modules in this decomposition which have 
dual
endpoint $t$ and diameter $d$. It is well-known that $mult(t,d)$ is
independent of the decomposition (cf. \cite{curtisReiner}).
\end{definition}

\begin{definition}\label{upsilondef}
With reference to Definition \ref{setup},
define a set $\Upsilon$ by
\begin{displaymath}
\Upsilon := \{ (t,d) \in \mathbb{Z}^2 \; | \; 0 \leq d
    \leq D,  \; \;  {\textstyle \frac{1}{2}}
        (D-d) \leq t \leq D-d \}.
\end{displaymath}
By (\ref{2t+deq}),
$\; mult(t,d)=0$ for all integers $t,d$ such that
$(t,d) \notin \Upsilon$.
We define a partial order $\preceq$ on $\Upsilon$ by
\begin{displaymath}
    (t,d) \preceq (t',d') \hspace{1em} \mbox{ if and only if } \hspace{1em}
            t \leq t' \; \mbox{ and } \; t'+d' \leq t+d.
\end{displaymath}
\end{definition}


\begin{example} 
    With reference to Definition \ref{setup}, suppose $D=7$. In
    Figure \ref{fig1}, we represent each element $(t,d) \in \Upsilon$
    by a line segment beginning in the $t^{th}$ column and having
    length $d$. For any elements $a \in \Upsilon$, $b \in \Upsilon$,
    observe that $a \preceq b$ if and only if  the line segment
    representing $a$ extends the line segment representing $b$.
\end{example}

\medskip

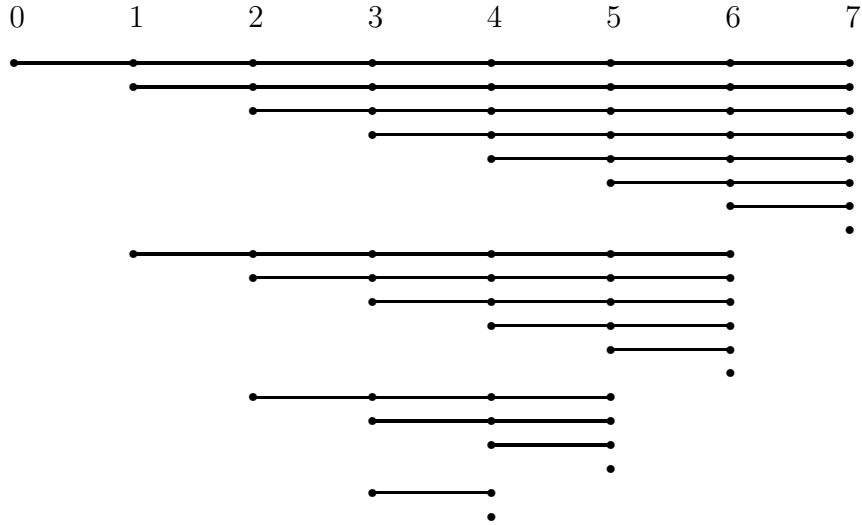
\begin{figur}{The set $\Upsilon$ when $D=7$}\label{fig1}

\bigskip
\bigskip

\center
\setlength{\unitlength}{0.012500in}%
\begingroup\makeatletter\ifx\SetFigFont\undefined
\def\x#1#2#3#4#5#6#7\relax{\def\x{#1#2#3#4#5#6}}%
\expandafter\x\fmtname xxxxxx\relax \def\y{splain}%
\ifx\x\y   
\gdef\SetFigFont#1#2#3{%
  \ifnum #1<17\tiny\else \ifnum #1<20\small\else
  \ifnum #1<24\normalsize\else \ifnum #1<29\large\else
  \ifnum #1<34\Large\else \ifnum #1<41\LARGE\else
     \huge\fi\fi\fi\fi\fi\fi
  \csname #3\endcsname}%
\else
\gdef\SetFigFont#1#2#3{\begingroup
  \count@#1\relax \ifnum 25<\count@\count@25\fi
  \def\x{\endgroup\@setsize\SetFigFont{#2pt}}%
  \expandafter\x
    \csname \romannumeral\the\count@ pt\expandafter\endcsname
    \csname @\romannumeral\the\count@ pt\endcsname
  \csname #3\endcsname}%
\fi
\fi\endgroup
\begin{picture}(354,202)(88,553)
\thicklines
\put(90, 755){\line( 1, 0){349}}
\put( 90,755){\circle*{3}}
\put(140,755){\circle*{3}}
\put(190,755){\circle*{3}}
\put(240,755){\circle*{3}}
\put(290,755){\circle*{3}}
\put(340,755){\circle*{3}}
\put(390,755){\circle*{3}}
\put(440,755){\circle*{3}}
\put(140,745){\line( 1, 0){299}}
\put(140,745){\circle*{3}}
\put(190,745){\circle*{3}}
\put(240,745){\circle*{3}}
\put(290,745){\circle*{3}}
\put(340,745){\circle*{3}}
\put(390,745){\circle*{3}}
\put(440,745){\circle*{3}}
\put(190,735){\line( 1, 0){249}}
\put(190,735){\circle*{3}}
\put(240,735){\circle*{3}}
\put(290,735){\circle*{3}}
\put(340,735){\circle*{3}}
\put(390,735){\circle*{3}}
\put(440,735){\circle*{3}}
\put(240,725){\line( 1, 0){199}}
\put(240,725){\circle*{3}}
\put(290,725){\circle*{3}}
\put(340,725){\circle*{3}}
\put(390,725){\circle*{3}}
\put(440,725){\circle*{3}}
\put(290,715){\line( 1, 0){149}}
\put(290,715){\circle*{3}}
\put(340,715){\circle*{3}}
\put(390,715){\circle*{3}}
\put(440,715){\circle*{3}}
\put(340,705){\line( 1, 0){99}}
\put(340,705){\circle*{3}}
\put(390,705){\circle*{3}}
\put(440,705){\circle*{3}}
\put(390,695){\line( 1, 0){49}}
\put(390,695){\circle*{3}}
\put(440,695){\circle*{3}}
\put(440, 685){\circle*{3}}
\put(140,675){\line( 1, 0){249}}
\put(140,675){\circle*{3}}
\put(190,675){\circle*{3}}
\put(240,675){\circle*{3}}
\put(290,675){\circle*{3}}
\put(340,675){\circle*{3}}
\put(390,675){\circle*{3}}
\put(190,665){\line( 1, 0){199}}
\put(190,665){\circle*{3}}
\put(240,665){\circle*{3}}
\put(290,665){\circle*{3}}
\put(340,665){\circle*{3}}
\put(390,665){\circle*{3}}
\put(240,655){\line( 1, 0){149}}
\put(240,655){\circle*{3}}
\put(290,655){\circle*{3}}
\put(340,655){\circle*{3}}
\put(390,655){\circle*{3}}
\put(290,645){\line( 1, 0){99}}
\put(290,645){\circle*{3}}
\put(340,645){\circle*{3}}
\put(390,645){\circle*{3}}
\put(340,635){\line( 1, 0){49}}
\put(340,635){\circle*{3}}
\put(390,635){\circle*{3}}
\put(390,625){\circle*{3}}
\put(190,615){\line( 1, 0){149}}
\put(190,615){\circle*{3}}
\put(240,615){\circle*{3}}
\put(290,615){\circle*{3}}
\put(340,615){\circle*{3}}
\put(240,605){\line( 1, 0){99}}
\put(240,605){\circle*{3}}
\put(290,605){\circle*{3}}
\put(340,605){\circle*{3}}
\put(290,595){\line( 1, 0){49}}
\put(290,595){\circle*{3}}
\put(340,595){\circle*{3}}
\put(340,585){\circle*{3}}
\put(240,575){\line( 1, 0){49}}
\put(240,575){\circle*{3}}
\put(290,575){\circle*{3}}
\put(290,565){\circle*{3}}
\put( 88,770){\makebox(0,0)[lb]{\smash{\SetFigFont{12}{14.4}{rm}0}}}
\put(138,770){\makebox(0,0)[lb]{\smash{\SetFigFont{12}{14.4}{rm}1}}}
\put(188,770){\makebox(0,0)[lb]{\smash{\SetFigFont{12}{14.4}{rm}2}}}
\put(238,770){\makebox(0,0)[lb]{\smash{\SetFigFont{12}{14.4}{rm}3}}}
\put(288,770){\makebox(0,0)[lb]{\smash{\SetFigFont{12}{14.4}{rm}4}}}
\put(338,770){\makebox(0,0)[lb]{\smash{\SetFigFont{12}{14.4}{rm}5}}}
\put(388,770){\makebox(0,0)[lb]{\smash{\SetFigFont{12}{14.4}{rm}6}}}
\put(438,770){\makebox(0,0)[lb]{\smash{\SetFigFont{12}{14.4}{rm}7}}}

\end{picture}

\end{figur}

\bigskip

\begin{lemma}\label{zz3}
With reference to Definition \ref{setup} and Definition 
\ref{upsilondef}, fix any $(t,d)
    \in \Upsilon$.  Then
\begin{equation}\label{ieiei3}
trace(E_{t} L^{*d} R^{*d} E_t) = m_{t} \prod_{h=t}^{t+d-1} b_h^{*} 
c_{t+d-h}^{*}.
\end{equation}
\end{lemma}

\begin{proof}
    By \cite[Lemma 4.1]{dickieNoteThin}, we find
\begin{equation}  \label{star}
    trace( E_tA^*_dE_{t+d}A^*_dE_t) = m_{t} \, q^{t}_{d,t+d}. 
\end{equation}   
By \cite[p. 276]{bannai} and since $\Gamma$ is Q-polynomial, we find
\begin{equation} \label{2star}
     q^{t}_{d,t+d} =
      \prod_{h=t}^{t+d-1} \frac{b^*_h}{c^*_{t+d-h}}.  
\end{equation}         
     We claim
\begin{equation}  \label{3star}
    E_{t+d}A^*_dE_t = \frac{1}{c_{1}^{*}c_{2}^{*}\cdots c_{d}^{*}}
    R^{*d}E_t.                
\end{equation}
     To see this, recall that for $0 \leq i \leq D$,
     $A^*_i$ is a polynomial in $A^*$ with degree $i$ and
     leading coefficient
     $(c^*_1c^*_2 \cdots c^*_i)^{-1}$.  Combining this and 
     (\ref{triangleQ}) we find
     $E_{t+d}A^{*i}E_t = 0$ for $0 \leq i \leq d-1$.
     We may now argue
\begin{eqnarray*}
     c^*_1c^*_2 \cdots c^*_d E_{t+d}A^*_dE_t &=&
            E_{t+d}A^{*d}E_t \\
             &=& E_{t+d}A^*E_{t+d-1}A^* \cdots E_{t+1}A^*E_t  \\
             &=& R^{*d}E_t.
\end{eqnarray*}
     We now have (\ref{3star}).  Taking the transpose of (\ref{3star}) we get
\begin{equation}  \label{4star}
     E_{t}A^*_dE_{t+d} = \frac{1}{c_{1}^{*}c_{2}^{*}\cdots c_{d}^{*}}
     E_tL^{*d}.
\end{equation}
Evaluating (\ref{star}) using
(\ref{2star})--(\ref{4star}), we get
the  desired result.
\end{proof}

\begin{lemma}
With reference to Definition \ref{setup} and Definition 
\ref{upsilondef},
fix $(i,j) \in \Upsilon$ and $(t,d) \in \Upsilon$.  Then
for any irreducible $T$-module $W$ with dual endpoint $i$ and diameter $j$,
\begin{equation}\label{ieei}
    trace(E_t L^{*d} R^{*d} E_{t}) \big|_W
            = \prod_{h=t-i}^{t-i+d-1} b_h^{*}(W) c_{h+1}^{*}(W)
\end{equation}
    if $(i,j) \preceq (t,d), \; $ and
        $ \; \; trace(E_t L^{*d} R^{*d} E_{t}) \big|_W =
    0 \; \; $ if $ \; \; (i,j) \not\preceq (t,d)$.
\end{lemma}

\begin{proof}
    Let $r$ denote the endpoint of $W$ and pick any nonzero $v \in E_r^{*} W$.
    By Theorem \ref{basis}(ii), $B:= \{ E_h v \; | \; \intv{i}h{i+j} \}$
    is a basis for $W$.  We consider the action of
    $E_t L^{*d} R^{*d} E_{t}$ on $B$.
    First assume $(i,j) \preceq (t,d)$. Then
    $E_t L^{*d} R^{*d} E_{t}$ vanishes on every element of $B$ except
    $E_tv$, and
    $$E_t L^{*d} R^{*d} E_{t} (E_t v) = \prod_{h=t-i}^{t-i+d-1}
                                    b_h^{*}(W) c_{h+1}^{*}(W) E_t v$$
    by (\ref{cswiparteq}), (\ref{bswiparteq}).
    Line (\ref{ieei}) follows.
    Next assume $(i,j) {\not \preceq} (t,d)$. Then
    $E_t L^{*d} R^{*d} E_{t}$ vanishes on each element of $B$,
    and so its trace on $W$ is zero.
\end{proof}

\medskip

    To state our next theorem, we need a bit of notation.
    Let $Y$ be as in Definition \ref{setup} and let $\Upsilon$ be as in
    Definition \ref{upsilondef}.
    Select $(t,d) \in \Upsilon$.
    We define $c_0^{*} (t,d) := 0, \; b_d^{*}(t,d) := 0$.
    For $d \geq 1$ we further define 
    \begin{eqnarray*}
    c_i^{*} (t,d) &:=& \frac{ ( \theta^2_{t+i} - \theta^2_t)
                        ( \theta^*_{r+2} - \theta^*_{r+1} ) +
                (\theta_t \theta_{t+1} - \theta_{t+i} \theta_{t+i+1} )
                        ( \theta^*_{r+1} - \theta^*_r ) }
    {(\theta_{t+i-1} - \theta_{t+i})(\theta_{t+i-1} -
    \theta_{t+i+1})}
                        \hspace{2em} (\interval1i{d-1}), \nonumber \\
    c^*_d (t,d) &:=& \frac{\theta_{t+d} (\theta^*_{r+1}-\theta^*_{r})}
                             {\theta_{t+d-1} - \theta_{t+d}},  \\
    \nonumber \\
    b^*_0 (t,d) &:=& \frac{\theta_t (\theta^*_r-\theta^*_{r+1})}
                             {\theta_t - \theta_{t+1}},  \\
    b_i^{*} (t,d) &:=&  \frac{ ( \theta^2_{t+i} - \theta^2_t)
                        ( \theta^*_{r+2} - \theta^*_{r+1} ) +
                (\theta_t \theta_{t+1} - \theta_{t+i} \theta_{t+i-1} )
                        ( \theta^*_{r+1} - \theta^*_r ) }
    {(\theta_{t+i+1} - \theta_{t+i})(\theta_{t+i+1} -
    \theta_{t+i-1})}
                        \hspace{2em} (\interval1i{d-1}), \nonumber \\
    \end{eqnarray*}
    where $r := D-d$.  Observe that if $W$ is any
    irreducible $T$-module with dual endpoint $t$ and diameter $d$,
    then $(t,d) \in \Upsilon$ and $c_i^{*}(t,d) = c_i^{*}(W)$,
    $b_i^{*}(t,d) = b_i^{*}(W)$ $(\intv0id)$.
    However, such a module need not exist.
    
    We now give a recurrence which will enable us to compute the 
    multiplicities of the irreducible $T$-modules.

\begin{theorem}\label{bigkahoona}
    With reference to Definition \ref{setup}, fix any $(t,d)
    \in \Upsilon$.  Then
\begin{equation}\label{ieei2}
 m_{t} \prod_{h=t}^{t+d-1} b_h^{*} 
c_{t+d-h}^{*}=  \sum_{(i,j) \in \Upsilon
                                            \atop (i,j) \preceq (t,d)}
                 mult(i,j) \prod_{h=t-i}^{t-i+d-1} b_h^{*}(i,j) c_{h+1}^{*}(i,j).
\end{equation}
\end{theorem}

\begin{proof}
Since $T$ is semisimple, we may
decompose the standard module $V$ as
$$
V = W_1 + W_2 + \cdots + W_s \hspace{3em} (\mbox{orthogonal direct sum}),
$$
where $W_1$, $W_2$,..., $W_s$ are irreducible $T$-modules.
It follows that
\begin{eqnarray}
trace(E_t L^{*d} R^{*d} E_{t}) &=& \sum_{i=1}^s trace(E_t L^{*d} R^{*d} E_{t})
\big{|}_{W_i} \label{7458fj}.
\end{eqnarray}
    Evaluating (\ref{7458fj}) using (\ref{ieiei3}), (\ref{ieei})
    we obtain (\ref{ieei2}).
\end{proof}

\begin{remark}\label{prudentrmk}
    With reference to Definition \ref{setup},
    we can use Theorem \ref{bigkahoona} to
    recursively compute the multiplicities
    $\{ mult(t,d) \; | \; (t,d) \in \Upsilon \}$.
    Indeed, pick any $(t,d) \in \Upsilon$. Then
    (\ref{ieei2}) gives a linear equation
    in the variables
    $\{ mult(i,j) \; | \; (i,j) \in \Upsilon, \; (i,j) \preceq (t,d) \}$.
    In this equation, the coefficient of $mult(t,d)$ is
\begin{equation}\label{donkeyfoot}
                 \prod_{h=0}^{d-1} b_h^{*}(t,d) c_{h+1}^{*}(t,d).
\end{equation}
    Suppose the coefficient (\ref{donkeyfoot}) is nonzero.  Then we can
    divide both sides of equation (\ref{ieei2})  by it, and obtain
    $mult(t,d)$ in terms of
    $\{ mult(i,j) \; | \; (i,j) \prec (t,d) \}$.
    Suppose the coefficient (\ref{donkeyfoot}) equals $0$.
    By Corollary \ref{8394}(ii), there is no module with dual endpoint $t$
    and diameter $d$, so $mult(t,d)=0$.
\end{remark}

\bigskip

We now illustrate Theorem \ref{bigkahoona} with an example.   

\begin{example}  \label{multlisting} With reference to Definition \ref{setup}, suppose
    $Y$ is not one of $O_{D+1}$, $\Box_{2D+1}$.  For all $(t,d) \in 
    \Upsilon$ such that $d \geq D-3$, the scalar $\; mult(t,d)$ is given 
    below. 
\begin{enumerate}
\begin{item}  $mult(0,D) =1$.
\end{item}

\begin{item}  $\displaystyle{mult(1,D-1) = 
\frac{(q^{2D}-1)(1+sq^{2})}{(1-q)(1+sq^{2D+1})}}$.
\end{item}

\begin{item}  $\displaystyle{mult(1, D-2) = 
\frac{(q^{2D}-q^{2})(1+sq)(1+sq^{2})(sq^{2D+2}-1)}{(q^{2}-1)(s^{2}q^{2D+4}-1)
(1+sq^{2D+1})}}$.
\end{item}

\begin{item} $\displaystyle{mult(2, D-2) = \frac{(q^{2D}-1)(q^{2D}-q^{2})
(1+sq)(1+sq^{4})(s^{2}q^{2D+3}-1)}{q(q+1)(q-1)^{2}(s^{2}q^{2D+4}-1)
(1+sq^{2D})(1+sq^{2D+1})}}$.
\end{item}

\begin{item} $\displaystyle{mult(2, D-3) = 
\frac{(q^{2D}-1)(q^{2D}-q^{4})(1+sq)(1+sq^2)(1+sq^4)
(1-sq^{2D+2})}{q(q-1)(q^{2}-1)(1+sq^{2D+1})(sq^{D+3}-1)(q+sq^{2D})(1+sq^{D+3})}}$.
\end{item}

\begin{item} $\displaystyle{mult(3, D-3) 
=\frac{(q^{2D}-1)(q^{2D}-q^{2})(q^{2D}-q^{4})
(1+sq)(1+sq^{2})(1+sq^{6})(1-s^{2}q^{2D+3})}{q^{2}(q-1)(q^{2}-1)(q^{3}-1)(1+sq^{D+3})
(sq^{D+3}-1)(q+sq^{2D})(1+sq^{2D})(1+sq^{2D+1})}}$.
\end{item}
\end{enumerate}
     The scalars $q,s$ are from Lemma \ref{qbipform}.  We remark the 
     denominators in (i)--(vi) are 
     nonzero by Corollary \ref{nonzeros} and Corollary \ref{bounds}.
\end{example}

\begin{proof}
     Solve (\ref{ieei2}) recursively for the multiplicities, as 
     outlined in Remark \ref{prudentrmk}.  Evaluate the results using 
     Lemmas \ref{qbipform}, 
\ref{qbipformII} and Corollaries \ref{newthi}, \ref{hhs}. 
\end{proof}

\end{document}